\documentclass[12pt,leqno,a4paper]{annacad}%
\pdfoutput=1
\usepackage{amssymb}
\usepackage{esint}
\usepackage{mathtools}
\usepackage{hyperref}
\usepackage{textcomp}
\mathtoolsset{showonlyrefs}
\usepackage[toc]{appendix}

\usepackage[latin,finnish,english]{babel}
\usepackage[latin1]{inputenc}
\usepackage[pdftex]{graphicx}

\usepackage{pgf,tikz}
\usetikzlibrary{arrows}

\theoremstyle{plain}
\newtheorem{theorem}{Theorem}

\newtheorem{question}[theorem]{Question}
\theoremstyle{definition}
\newtheorem{definition}[theorem]{Definition}

\theoremstyle{remark}

\renewcommand{\contentsname}{Contents of the introductory part}

\title{On the broken ray transform}
\author{Joonas Ilmavirta}

\usepackage{pdfpages}

\begin{document}


\includepdf[pages={3-4}]{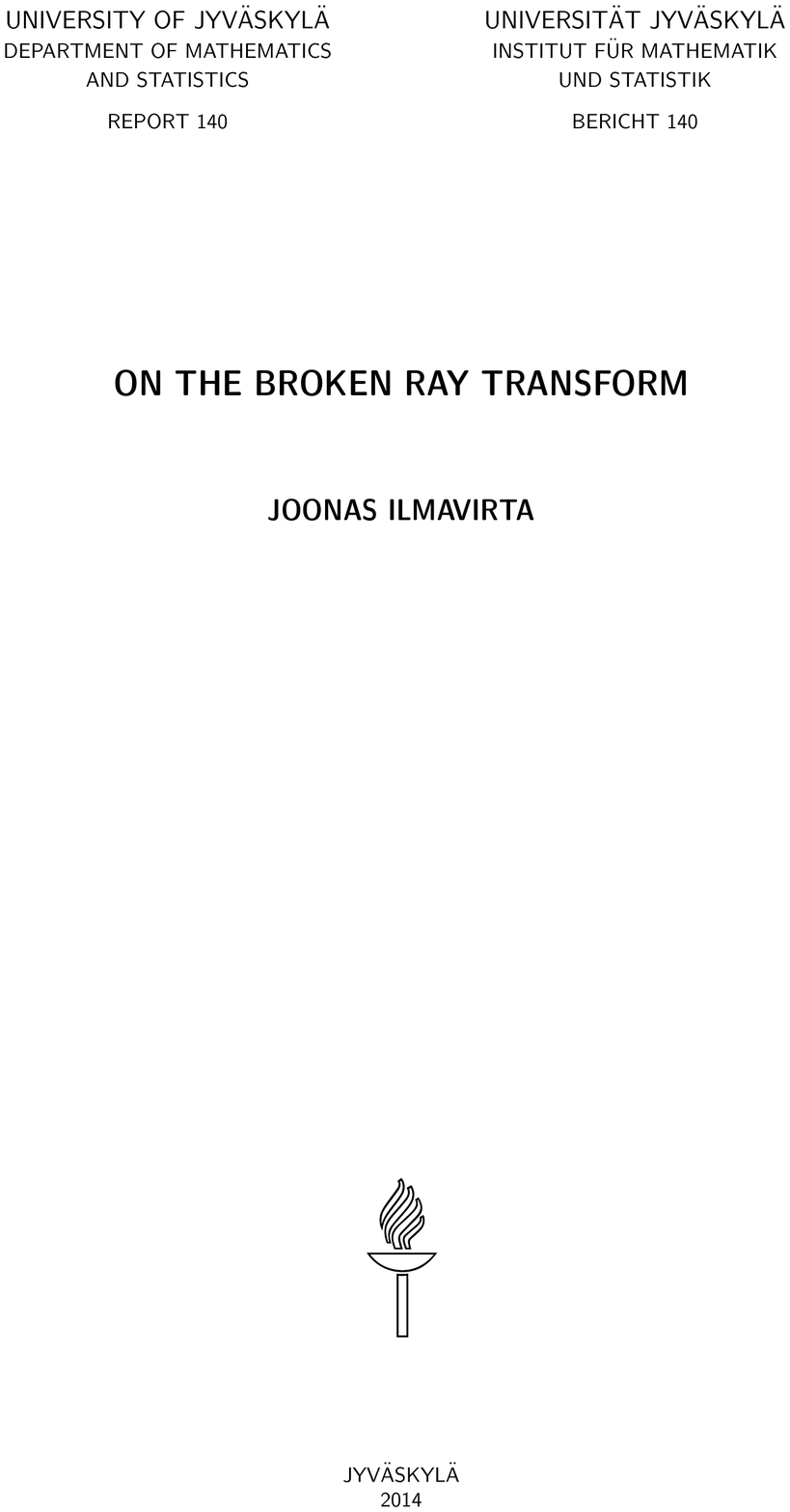}

\newpage

\section*[Abstract]{Abstract}

The classical problem of X-ray tomography asks whether one can reconstruct a function from its integrals over all lines.
This problem has a wide range of applications and has been studied intensively.
As an example of the applications we mention computerized tomography, which is an important tool in medical imaging.
When we replace the lines with broken rays which reflect on some subset of the boundary of the domain of interest, the problem becomes much more difficult.
The reflecting part of the boundary is inaccessible to measurements, so the problem at hand resembles X-ray tomography with partial data.
This type of tomography is called broken ray tomography in this thesis.

The integral transform corresponding to broken ray tomography is called the broken ray transform, which can be thought of as a generalization of the X-ray transform.
The fundamental question is whether this transform is injective.
We employ four different methods to approach this question, and each of them gives interesting results.

Direct calculation can be used in a ball, where the geometry is particularly simple.
If the reflecting part of the boundary is (piecewise) flat, a reflection argument can be used to reduce the problem to the usual X-ray transform.
In some geometries one can use broken rays near the boundary to determine the values of the unknown function at the reflector, and even construct its Taylor series.
One can also use energy estimates -- which in this context are known as Pestov identities -- to show injectivity in the presence of one convex reflecting obstacle.
Many of these methods work also on Riemannian manifolds.

We also discuss the periodic broken ray transform, where the integrals are taken over periodic broken rays.
The broken ray transform and its periodic version have applications in other inverse problems, including Calder\'on's problem and problems related to spectral geometry.

A simplified introduction to inverse problems to the nonmathematical reader is given in appendices
\ref{sec:eng}~(English),
\ref{sec:fin}~(Finnish), and
\ref{sec:lat}~(Latin).

\newpage

\begin{otherlanguage}{finnish}

\section*[Tiivistelm‰]{Tiivistelm‰}

Rˆntgen-tomografian klassinen ongelma on selvitt‰‰, voiko funktion rekonstruoida sen integraaleista kaikkien suorien yli.
T‰ll‰ ongelmalla on runsaasti sovelluksia ja sit‰ on tutkittu ahkerasti.
Sovelluksista mainittakoon tietokonetomografia, joka on t‰rke‰ tyˆkalu l‰‰ketieteellisess‰ kuvantamisessa.
Kun suorat korvataan kiinnostavan alueen reunan jostain osasta heijastuvilla murtos‰teill‰, ongelma muuttuu paljon vaikeammaksi.
Reunan heijastava osa ei ole k‰ytett‰viss‰ mittauksiin, joten ongelma muistuttaa osittaisen datan Rˆntgen-tomografiaa.
T‰t‰ tomografiatyyppi‰ kutsutaan t‰ss‰ tyˆss‰ murtos‰detomografiaksi.

Murtos‰detomografiaa vastaava integraalimuunnos on murtos‰demuunnos, jota voi ajatella Rˆntgen-muunnoksen yleistyksen‰.
Peruskysymys on, onko t‰m‰ muunnos injektiivinen.
K‰yt‰mme nelj‰‰ eri tapaa l‰hesty‰ t‰t‰ ongelmaa, ja niist‰ jokainen antaa mielenkiintoisia tuloksia.

Suoraa laskua voi k‰ytt‰‰ pallossa, jossa geometria on poikkeuksellisen yksinkertainen.
Jos reunan heijastava osa on (paloittain) laakea, voi ongelman heijastusargumentilla palauttaa tavalliseen Rˆntgen-muunnokseen.
Joissain geometrioissa voi k‰ytt‰‰ l‰hell‰ reunaa olevia murtos‰teit‰ m‰‰ritt‰m‰‰n tuntemattoman funktion heijastavalla osalla ja jopa konstruoimaan sen Taylorin sarjan.
On myˆs mahdollista k‰ytt‰‰ energiaestimaatteja -- jollaisia t‰ss‰ yhteydess‰ kutsutaan Pestovin identiteeteiksi -- injektiivisyyden todistamisessa, kun alueessa on yksi konveksi heijastava este.
Monet n‰ist‰ menetelmist‰ toimivat myˆs Riemannin monistoilla.

Tarkastelemme myˆs jaksollista murtos‰demuunnosta, jossa integrointi tehd‰‰n kaikkien jaksollisten murtos‰teiden yli.
Murtos‰demuunnoksella ja sen jaksollisella versiolla on sovelluksia muissa inversio-ongelmissa kuten Calder\'onin ongelmassa ja spektraaligeometriaan liittyviss‰ kysymyksiss‰.

Yksinkertaistettu johdanto inversio-ongelmiin matematiikkaa tuntemattomalle lukijalle on liitteiss‰
\ref{sec:eng}~(englanniksi),
\ref{sec:fin}~(suomeksi) ja
\ref{sec:lat}~(latinaksi).

\end{otherlanguage}

\newpage

\section*[Acknowledgements]{Acknowledgements}

I wish to express my greatest gratitude to my advisor, Professor Mikko Salo, for leading me to an interesting path and for giving support and advice on my journey along it.
I have been able to complete my studies in a friendly and supportive atmosphere, for which I thank all of the staff and students at the Department of Mathematics and Statistics.
Numerous colleagues have helped me with the problems I encountered, and I hope to have been able to offer them something in return.

I am most thankful to Professor Vladimir Sharafutdinov for making invaluable remarks on an early version of this thesis and for recommending several relevant articles that I was previously unaware of.

I would never have finished this work without life outside academia.
Support from family and friends has made all of this possible.
The people who have been there for me know who they are; I refrain from listing names as it would only lead to me having to decide whom not to include in the list.

My research was financially supported by the Academy of Finland, especially the Centre of Excellence in Inverse Problems Research, for which I am most grateful.

\vspace{2em}
\begin{flushright}
Jyv\"askyl\"a, March 2014\\
Joonas Ilmavirta
\end{flushright}

\newpage

\section*[List of included articles]{List of included articles}

\noindent
This thesis is based on the work contained within the following publications:

%
%
%
%
%
%

\vspace{1em}
\begingroup
\renewcommand{\section}[2]{}

\endgroup

\noindent
The author has participated actively in development of the joint paper~\cite{IS:brt-pde-1obst}.
In particular, the author is responsible for the regularity lemma and calculating the linearization.

The Latin introduction (appendix~\ref{sec:lat}) has been published in Melissa (number 180, June 2014), a magazine published in Latin.

\newpage

\renewcommand*{\contentsname}{Contents of the introductory part}

\setcounter{tocdepth}{1}
\tableofcontents

\noindent
Appendices~\ref{sec:eng}, \ref{sec:fin} and~\ref{sec:lat} contain an introduction to inverse problems for the nonmathematical audience.
The same content is presented in three languages.

\newpage

\newcommand{\R}{\mathbb R}
\newcommand{\C}{\mathbb C}
\newcommand{\N}{\mathbb N}
\newcommand{\Z}{\mathbb Z}
\newcommand{\Q}{\mathbb Q}
\newcommand{\h}{\mathcal H}
\newcommand{\A}{\mathcal A}
\newcommand{\At}{\widetilde{\mathcal A}}
\newcommand{\brt}{{\mathcal G}}
\newcommand{\rt}{{\mathcal I}}
\newcommand{\q}{\mathcal Q}
\newcommand{\F}{\mathcal F}
\newcommand{\der}{\mathrm{d}}
\newcommand{\eps}{\varepsilon}
\renewcommand{\phi}{\varphi}
\newcommand{\abs}[1]{\left| #1 \right|}
\newcommand{\aabs}[1]{\left\| #1 \right\|}
\newcommand{\sisus}{\operatorname{int}}
\newcommand{\spt}{\operatorname{spt}}
\newcommand{\co}{\operatorname{co}}
\newcommand{\tr}{\operatorname{tr}}
\newcommand{\Lip}{\operatorname{Lip}}
\newcommand{\lip}{\operatorname{lip}}
\newcommand{\dive}{\operatorname{div}}
\renewcommand{\theta}{\vartheta}
\newcommand{\ang}{{\mathrm{ang}}}
\newcommand{\syt}{{\mathrm{syt}}}
\newcommand{\qa}{{\mathrm{qa}}}
\newcommand{\rot}{\mathfrak R}
\newcommand{\radon}{\mathcal R}
\newcommand{\braket}[2]{\left\langle #1 \middle| #2 \right \rangle}
\newcommand{\Order}{\mathcal O}
\newcommand{\order}{o}
\newcommand{\rfl}[1]{\widetilde{#1}}
\newcommand{\ext}[1]{\widehat{#1}}
\newcommand{\id}{\operatorname{id}}
\newcommand{\fun}[2]{{#2}^{#1}}
\newcommand{\lap}{\mathcal L}
\newcommand{\sff}[2]{\mathrm{I\!I}\!\left(#1,#2\right)}
\newcommand{\tpm}[1]{\sff{#1}{#1}}
\newcommand{\ip}[2]{\left\langle#1,#2\right\rangle}
\newcommand{\iip}[2]{\left(#1,#2\right)}
\newcommand{\bd}[1]{{\check{#1}}}
\newcommand{\dbd}[1]{\dot{\bd{#1}}}
\newcommand{\ddbd}[1]{\ddot{\bd{#1}}}
\newcommand{\Der}[1]{\frac{\der}{\der #1}}
\newcommand{\roo}[1]{{\tilde{#1}}}
\newcommand{\lin}{\mathcal L}

\section{Introduction}
\label{sec:intro}

\subsection{Inverse problems}
\label{sec:ip}

To understand an inverse problem, we first need to understand direct problems.
Consider, for example, the following direct problems:
\begin{itemize}
\item Given a~$C_0$ function (continuous with compact support) $f:\R^n\to\R$, find its integral over every line.
\item Given a closed Riemannian manifold, find the spectrum of its Laplace-Beltrami operator.
\item Given a compact Riemannian manifold with boundary, find the distance between any two boundary points.
\item Given a~$C^2$ function $\gamma:\R^n\supset\bar\Omega\to(0,\infty)$, find the boundary values $u|_{\partial\Omega}$ and $\gamma\partial_\nu u|_{\partial\Omega}$ for all solutions $u\in W^{1,2}(\Omega)$ of $\dive(\gamma\nabla u)=0$.
\end{itemize}
It is obvious that the solution to these problems exists uniquely, and in some cases it is relatively easy to find an explicit solution.

When turned around, these become inverse problems:
\begin{itemize}
\item Given the integral of a~$C_0$ function over every line, find the function. (The X-ray tomography problem.)
\item Given the spectrum of the Laplace-Beltrami operator on a closed manifold, find the manifold. (Related to the spectral rigidity problem.)
\item Given the distance between any two boundary points on a manifold with boundary, find the manifold. (The boundary rigidity problem.)
\item Given the boundary values $u|_{\partial\Omega}$ and $\gamma\partial_\nu u|_{\partial\Omega}$ for all solutions $u\in W^{1,2}(\Omega)$ of $\dive(\gamma\nabla u)=0$, find the function $\gamma\in C^2(\bar\Omega)$. (Calder\'on's problem.)
\end{itemize}
These famous inverse problems have been studied intensively.
Solving these problems is far more difficult than solving the corresponding direct problems.
In particular, it is not clear that a solution should exist, let alone be unique.
If the solution exists uniquely, one would like to have an explicit reconstruction method and know that it is stable.


\subsection{Ray transforms}
\label{sec:rt}

We are interested in inverse problems similar to the first one in the list above.
The question whether on can recover a function in a bounded Euclidean domain from its integral over all lines has applications in medical imaging.
This problem was first solved by Radon~\cite{radon} and Cormack~\cite{cormack}.
For applications and further developments, see for example the book by Natterer~\cite{book-natterer}.

This question can be generalized in several directions.
First, lines can be replaced with another family of curves.
This includes geodesics on a Riemannian manifold.
Second, the integrals can be taken over hyperplanes instead of lines (the Radon transform) or more generally over all $d$-planes for any integer~$d$ strictly between one and the dimension (the $d$-plane Radon transform).
This second generalization does not make sense on all manifolds, but it can be defined in the Euclidean space and flat tori; the $d$-plane Radon transform on tori is discussed in~\cite{I:torus}.
Third, the scalar function may be replaced with a tensor field (see~\cite{S:tensor-book}).

The first generalization is the one most important to us.
Namely, we may replace the lines (or geodesics) with broken rays which are piecewise linear (or geodesic) curves.
There are many different types of breaking points to consider.
If the breaking points are in the interior of the domain, one can consider geodesics that break in all possible directions~\cite{KLU:broken-geodesic-flow} or only those that break in a fixed angle~\cite{MNTZ:radon,TN:radon,A:radon,AM:v-line-disc}.
Both of these cases have applications in imaging.
Recently Zhao, Schotland and Markel~\cite{ZSM:star-transform} (see also~\cite{KK:brt-inversion-range}) have also considered the case where a line breaks into several lines, forming a star rather than a broken line.
One may also let the broken rays reflect from the boundary of the domain according to the usual law of geometrical optics.
Versions of this problem have been considered by Mukhometov~\cite{M1,M2,M3,M4,M5,M6} and Eskin~\cite{eskin}, and this case is the main topic of this thesis.

%

\subsection{The broken ray transform}
\label{sec:brt}

Before going further, we need to introduce some notation.
Let~$M$ be a compact Riemannian manifold with boundary, and let its boundary be divided into disjoint sets~$E$ and~$R$.
A broken ray is a curve on~$M$ which is geodesic in the interior of~$M$, has both endpoints on~$E$ and reflects on other boundary points according to the usual reflection law.
That is, the incoming direction~$v_i$ and the outgoing direction~$v_o$ are related via $v_o=v_i-2\ip{\nu}{v_i}\nu$, where~$\nu$ is the outer unit normal at the reflection point.
In dimension two this amounts to demanding that the angle of incidence equals the angle of reflection.
Note that a broken ray \emph{must} reflect when it hits~$R$, but on~$E$ it can either reflect or stop.
Allowing reflections on~$E$ will be convenient.

We refer to the set~$E$ as the \emph{set of tomography}, since measurements can only be done on it.
The set~$R$ is a reflecting part of the boundary which is inaccessible to measurements.

There are two main types of geometrical situations that we are interested in.
The boundary~$\partial M$ can be connected, when~$E$ and~$R$ necessarily meet at some boundary point.
If the boundary is not connected,~$E$ and~$R$ may consist of connected components of the boundary, making them well separated.

For the sake of simplicity, we assume all broken rays to have finite length.
We allow reflections to be tangential (in which case the direction is not changed at all), although it will pose a problem that the reflection fails to be smooth at tangential directions.

\begin{definition}
\label{def:brt}
Let~$M$ be a manifold as described above and let~$\Gamma$ be the set of all broken rays on it.
Denote by~$\R^\Gamma$ the set of all mappings~$\Gamma\to\R$.
We define the \emph{broken ray transform} as a mapping $\brt:C(M)\to\R^\Gamma$ by letting
\begin{equation}
\brt f(\gamma)
=
\int_0^L f(\gamma(t))\der t
\eqqcolon
\int_{\gamma}f\der s
\end{equation}
for a broken ray $\gamma:[0,L]\to M$.
\end{definition}

The transform can be also be defined for many noncontinuous functions in the same way.
One can also introduce a weight or attenuation to the broken ray transform (see~\cite[Section~2.1]{I:bdy-det}), but for the sake of simplicity we omit them here.

\begin{question}
\label{q:brt}
Regarding the broken ray transform defined in definition~\ref{def:brt} above, we ask the following questions:
\begin{itemize}
\item Uniqueness: Is a function uniquely determined by its broken ray transform? In other words, is the broken ray transform injective?
\item Reconstruction: Is there an explicit numerical algorithm or analytic formula to recover a function from its broken ray transform?
\item Stability: If two functions have almost the same broken ray transform, are the functions almost the same?
\item Regularity: How do the answers to these questions depend on the regularity of the unknown function and the domain?
\end{itemize}
\end{question}

The goal of this thesis is to answer these questions.
A complete answer is by far beyond reach, but partial answers have been obtained.
We will give several examples of domains where the broken ray transform is injective.
We use minimal regularity assumptions for the methods of proof used, but there is no indication that the obtained regularity is optimal.
A partial stability result is shown in one case (see section~\ref{sec:pde}), but otherwise stability and reconstruction are only inherited by reducing the broken ray transform to the usual X-ray transform.

We also give counterexamples to uniqueness.
But since none of the counterexamples is particularly close to the positive results, we cannot claim that any of the theorems would be sharp.

\subsection{The periodic broken ray transform}
\label{sec:pbrt}

As the name suggests, a periodic broken ray is a periodic curve on a manifold~$M$ with boundary which reflects at the boundary according to the aforementioned reflection law.
Note that in this case all of the boundary~$\partial M$ is reflective, since periodic broken rays have no endpoints.

The periodic broken ray transform is then defined analogously with the broken ray transform.
We pose question~\ref{q:brt} also for the periodic broken ray transform.
The answers given here are limited to some examples and counterexamples.

\section{How to approach the broken ray transform}
\label{sec:howto}

We present four methods to approach the broken ray transform described in section~\ref{sec:brt}:

\begin{enumerate}
\item Compute the broken ray transform explicitly.
\item Reduce to X-ray transform by reflection.
\item Reduce boundary determination for the broken ray transform to injectivity of the X-ray transform on boundary via broken rays staying near the boundary.
\item Reduce injectivity of the broken ray transform to the unique solvability of a PDE.
\end{enumerate}

Details of these methods are given in the correspondingly numbered subsections~\ref{sec:explicit}--\ref{sec:pde}, but we give a short overview here.
The proofs given in the subsections are not complete, but merely attempt to convey the main ideas.
For detailed proofs, see the corresponding articles included in this thesis.

The first approach requires the domain to have very explicit geometry, at least for the reflecting part~$R$ of the boundary.
Such approach can be taken in a Euclidean disc, where a broken ray consists of rotated copies of a chord.
This approach was taken in~\cite{I:disk} and will be outlined in section~\ref{sec:explicit}.

The second approach is best illuminated with the example of a half plane:
If $f:\R\times[0,\infty)\to\R$ is a compactly supported continuous function, it can be reflected to a function $\tilde{f}:\R^2\to\R$ by setting $\tilde{f}(x,y)=f(x,\abs{y})$.
If we know the integral of~$f$ over all broken rays that reflect on $\R\times\{0\}$, then we know the integral of~$\tilde{f}$ over all lines in~$\R^2$.
Thus the broken ray transform is reduced to the X-ray transform, for which injectivity is known.

This approach has been taken in~\cite{I:disk,H:square,H:brt-flat-refl}.
The reflection method works best when the reflecting part~$R$ of~$\partial\Omega$ is flat.
More details of this method (including the role of flatness) are given in section~\ref{sec:refl}.

The third approach is based on an intuitive geometrical observation:
If~$\partial\Omega$ is strictly convex, then a broken ray starting almost tangent to~$\partial\Omega$ remains almost tangent to~$\partial\Omega$, and if a sequence of broken rays becomes closer and closer to being tangent to the boundary, the limit curve is a geodesic on~$\partial\Omega$.
In fact, any boundary geodesic along which~$\partial\Omega$ is strictly convex can be approached uniformly with broken rays.

Thus, if the unknown function is continuous, we may recover its integral over all boundary geodesics from the broken ray data.
If we can invert the X-ray transform on the manifold $\partial\Omega\setminus E$, we can reconstruct the unknown function at~$\partial\Omega$.
A similar argument can be used to recover the normal derivatives at the boundary as well, but the corresponding X-ray transform on the boundary manifold picks a weight depending on curvature.
A more detailed account of this method is given in section~\ref{sec:bdy-det}, which is based on~\cite{I:bdy-det}.

The last approach is not as obvious as the three other ones.
For an unknown function $f:\Omega\to\R$ we define another function~$u^f$ on the sphere bundle~$S\Omega$.
Assuming that~$f$ has a vanishing broken ray transform, this function~$u^f$ is shown to satisfy a PDE with some boundary conditions on~$E$ and~$R$.
An energy estimate (called Pestov identity in this context) can be used to show that this PDE has a unique solution and thus~$u^f$ vanishes.
The function~$u^f$ is defined so that this implies that~$f$ vanishes.

This kind of a PDE approach has been previously used for the X-ray transform (see for example~\cite{M:bdy-rig-surface-eng,PSU:tensor-survey}, the first of which is translated from~\cite{M:bdy-rig-surface-rus}) and also for the broken ray transform under very special geometrical assumptions by Eskin~\cite{eskin}.
Mukhometov~\cite{M1,M2,M3,M4,M5,M6} has also used this approach to study the broken ray transform.
This approach was used in~\cite{IS:brt-pde-1obst} to show injectivity of the broken ray transform on some Riemannian surfaces with one convex obstacle.
This approach can also be used to obtain a stability estimate: if the broken ray transform of~$f$ is small, then~$f$ itself is small.
All the necessary definitions and basic properties related to this method are given in section~\ref{sec:pde}.

The main result in~\cite{IS:brt-pde-1obst} is very close to that of~\cite{M4}.
A major motivation for~\cite{IS:brt-pde-1obst} is providing a basis for analysis of the broken ray transform in the presence of several convex obstacles.
Therefore regularity of the function~$u^f$ is studied in detail, to which end we also introduce Jacobi fields along broken rays.

It is important to note that the different approaches work in different geometrical settings.
The reflection approach works best when the reflecting part~$R$ of~$\partial\Omega$ is flat; the Riemannian metric of the constructed auxiliary manifold will fail to be~$C^1$ precisely when the second fundamental form does not vanish identically on~$R$.
The boundary determination result presented in section~\ref{sec:bdy-det} requires that~$R$ is strictly convex along sufficiently many boundary geodesics -- preferably all of them.
The PDE approach requires that~$R$ is strictly concave, and the reason is twofold: the billiard map is more regular in this case and an important term in the Pestov identity~\eqref{eq:pestov} has wrong sign otherwise.

We have thus ways to approach the broken ray transform if~$R$ is strictly convex, flat, or strictly concave.
When different parts of~$R$ have types, the methods may be combined to some extent: the reflection argument can easily be combined with any other one.
These reflector geometries do not, however, exhaust all cases.
Saddle-like~$R$, for example, does not yield to any of these approaches.

\subsection{Explicit calculation}
\label{sec:explicit}

The case when~$\Omega$ is a Euclidean disc is geometrically very convenient.
We may take~$\Omega$ to be the unit disc centered at the origin.
A broken ray in the disc consists of rotated copies of a chord.
Also, we may use polar coordinates $(r,\theta)$ and write functions as Fourier series in the angular variable~$\theta$.

This section follows the article~\cite{I:disk}.
The important property of the disc used here is rotational symmetry in polar coordinates.

We begin by observing that the broken ray transform in an $n$-dimensional Euclidean ball~$B^n$ for $n\geq2$ can be reduced to the case $n=2$.
If~$F$ is a two-dimensional subspace of~$\R^n$, $n\geq3$, then any broken ray with one segment in~$F$ is completely contained in $F$.
If we know injectivity for $n=2$, we know that the broken ray transform of $f:B^n\to\R$ determines~$f$ in $F\cap B^n$.
And since this holds for any~$F$, we have injectivity for every $n\geq2$.
We will therefore only prove the following theorems in the case $n=2$.

The main results obtained in~\cite{I:disk} are the following:

\begin{theorem}
\label{thm:disk-one}
If $f:\bar B^n\to\C$ is continuous and~$E$ is a singleton, the broken ray transform of~$f$ uniquely determines the integral of~$f$ over any circle centered at the origin with the singleton in the circle's plane.
In particular, if $n=2$, the integral of~$f$ is determined over any circle centered at the origin.
\end{theorem}

\begin{theorem}
\label{thm:disk-open}
Suppose $f:\bar B^n\to\C$ is uniformly quasianalytic in the angular variable in the sense of~\cite[Definition~14]{I:disk} and its angular derivatives of all orders satisfy the Dini-Lipschitz condition.
If~$E$ is open,~$f$ is uniquely determined by its broken ray transform.
\end{theorem}

A simple example of a function in two dimensions satisfying the regularity assumption of theorem~\ref{thm:disk-open} is
\begin{equation}
f(r,\theta)=\sum_{\abs{k}\leq K}a_k(r)e^{ik\theta},
\end{equation}
where $K\in\N$ and each~$a_k$ is H\"older continuous.

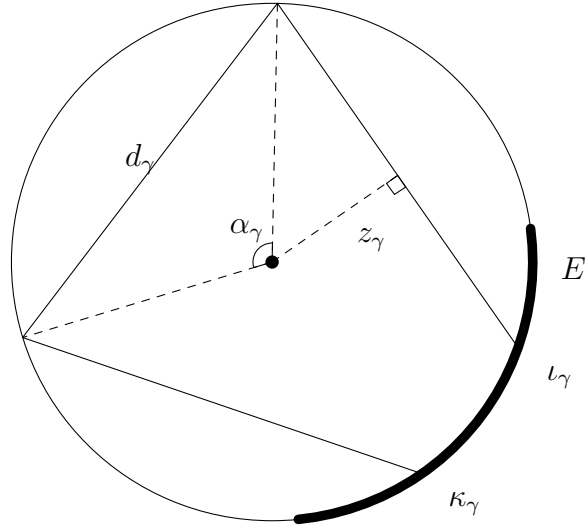
\begin{figure}%
\begin{tikzpicture}[line cap=round,line join=round,>=triangle 45,x=0.7cm,y=0.7cm]
\clip(-10.55,-0.11) rectangle (0.48,10.1);
\draw [shift={(-5.39,5.01)}] (0,0) -- (88.78:0.36) arc (88.78:196.93:0.36) -- cycle;
\draw(-3.24,6.5) -- (-3.09,6.29) -- (-2.88,6.43) -- (-3.03,6.65) -- cycle; 
\draw (-6.39,6.02) node[anchor=north west] {$ \alpha_\gamma $};
\draw(-5.39,5.01) circle (4.9);
\draw [dash pattern=on 3pt off 3pt] (-5.39,5.01)-- (-3.03,6.65);
\draw (-0.76,3.38)-- (-5.29,9.91);
\draw (-5.29,9.91)-- (-10.08,3.58);
\draw (-10.08,3.58)-- (-2.57,1);
\draw [dash pattern=on 3pt off 3pt] (-5.29,9.91)-- (-5.39,5.01);
\draw [dash pattern=on 3pt off 3pt] (-5.39,5.01)-- (-10.08,3.58);
\draw [shift={(-5.39,5.01)},line width=3.6pt]  plot[domain=-1.47:0.13,variable=\t]({1*4.9*cos(\t r)+0*4.9*sin(\t r)},{0*4.9*cos(\t r)+1*4.9*sin(\t r)});
\draw (-4,5.86) node[anchor=north west] {$ z_\gamma $};
\draw (-0.43,3.25) node[anchor=north west] {$ \iota_\gamma $};
\draw (-2.28,0.86) node[anchor=north west] {$ \kappa_\gamma $};
\draw (-8.36,7.46) node[anchor=north west] {$ d_\gamma $};
\draw (-0.17,5.31) node[anchor=north west] {$ E $};
\begin{scriptsize}
\fill [color=black] (-5.39,5.01) circle (2.5pt);
\end{scriptsize}
\end{tikzpicture}
\caption{An example of a broken ray~$\gamma$ in the disc. Here $m_\gamma=1$ and $n_\gamma=3$. The parameters~$\alpha_\gamma$, $z_\gamma$, $d_\gamma$, $\iota_\gamma$, and~$\kappa_\gamma$ as well as the set of tomography~$E$ are illustrated in the figure. The angles~$\iota_\gamma$ and~$\kappa_\gamma$ are identified with the corresponding points on the boundary.~\cite{I:disk}}%
\label{fig:br-param}%
\end{figure}

Both of these results are based on passing to the limit of infinitely many reflections.
To explain this in more detail, we must fix some notation.
We identify angles with points on the boundary.
A broken ray~$\gamma$
\begin{itemize}
\item contains~$n_\gamma$ line segments of length~$d_\gamma$ and distance~$z_\gamma$ from the origin,
\item winds around the origin~$m_\gamma$ times,
\item has initial and final points~$\iota_\gamma$ and~$\kappa_\gamma$, and
\item each line segment opens at an angle~$\alpha_\gamma$ as viewed from the origin.
\end{itemize}
These seven parameters define~$\gamma$ uniquely, and although this parametrization is redundant, it is intuitive and useful.
The winding number~$m_\gamma$ is defined by the relation
\begin{equation}
\label{eq:winding}
n_\gamma\alpha_\gamma
=
2\pi m_\gamma+\kappa_\gamma-\iota_\gamma.
\end{equation}
The lengths can be expressed in terms of the other parameters: $z_\gamma=\cos(\alpha_\gamma/2)$ and $d_\gamma=2\abs{\sin(\alpha_\gamma/2)}$. 
The parameters are illustrated in figure~\ref{fig:br-param}.

Suppose we have a sequence~$(\gamma_i)$ of broken rays such that $z_{\gamma_i}\to z\in(0,1)$ and $n_{\gamma_i}\to\infty$ as $i\to\infty$.
Broken rays with many reflections are almost spherically symmetric: the integral of~$f$ over~$\gamma_i$ tends to the integral of~$f$ against a spherically symmetric function.
For theorem~\ref{thm:disk-one} this observation is enough, since it turns the broken ray transform $f:B^2\to\C$ essentially to the Abel transform of the spherical mean function $a_0(r)=\fint_{\abs{x}=r}f(x)$.
Injectivity of the Abel transform finishes the proof.

Theorem~\ref{thm:disk-open} requires more.
First we write~$f$ as a Fourier series:
\begin{equation}
\label{eq:fourier}
f(r,\theta)=\sum_{k\in\Z}e^{ik\theta}a_k(r).
\end{equation}
Under the assumption of Dini-Lipschitz-continuity this series converges uniformly and therefore the broken ray transform may be computed term by term.

Another important observation is based on rotational symmetry:
If~$E$ is open, a slightly rotated broken ray is still a broken ray with endpoints in the same set~$E$.
But since rotating a broken ray clockwise has precisely the same effect on the broken ray transform than rotating the function counterclockwise, small rotations retain vanishing broken ray transform for a slightly smaller set of tomography.
Passing to the limit of zero rotation we find that the angular derivative of a function with vanishing broken ray transform has vanishing broken ray transform.
The argument can be iterated to show that angular derivatives of all orders have vanishing broken ray transform.

In theorem~\ref{thm:disk-one} we only recovered the Fourier component~$a_0$ by using the Abel transform (denoted here by~$\A_0$).
In attempt to recover all components and thus the whole function, we end up with generalized Abel transforms~$\A_k$.
These transforms are also injective on continuous functions and arise naturally in the Radon transform: if $f(r,\theta)=e^{ik\theta}a(r)$, then the integral of~$f$ over the line $\{x;x_1\cos\phi+x_2\sin\phi=\rho\}$ is $e^{ik\phi}\A_{\abs{k}}a(\rho)$ if $\rho>0$ and~$a$ is continuous.

\subsection{Reflection}
\label{sec:refl}

The only important difference between the X-ray transform and the broken ray transform is the existence of reflections.
As reflections can be difficult to handle in generic geometry, it is lucrative to look for methods that get rid of them altogether.

We reflect the domain~$\Omega$ at~$R$ to construct a new domain~$\tilde\Omega$ so that the natural projection map $\tilde\Omega\to\Omega$ takes lines into broken rays.
This construction is most clear in the half plane example given in section~\ref{sec:howto} where $\Omega=\R\times[0,\infty)$ and $\tilde\Omega=\R^2$.
This correspondence between broken rays in~$\Omega$ and lines in~$\tilde\Omega$ allows one to reduce injectivity (and reconstruction and stability) of the broken ray transform on~$\Omega$ to that of the X-ray transform on~$\tilde\Omega$.
If $R$ is something else than part of a hyperplane in a Euclidean space, the reflected domain~$\Omega$ is not Euclidean.
We recognize two types of ``non-hyperplane geometry'' for the reflecting part:~$R$ may be curved or may contain corners.

The reflection construction may also be carried out when~$\Omega$ is a Riemannian manifold -- this is done in the article~\cite{I:refl}.
We wish the new manifold~$\tilde\Omega$ to be such a manifold that the X-ray transform on it is injective.
Known injectivity results for Riemannian manifolds assume that the metric and the boundary are smooth.

\begin{figure}%
\includegraphics[scale=.5]{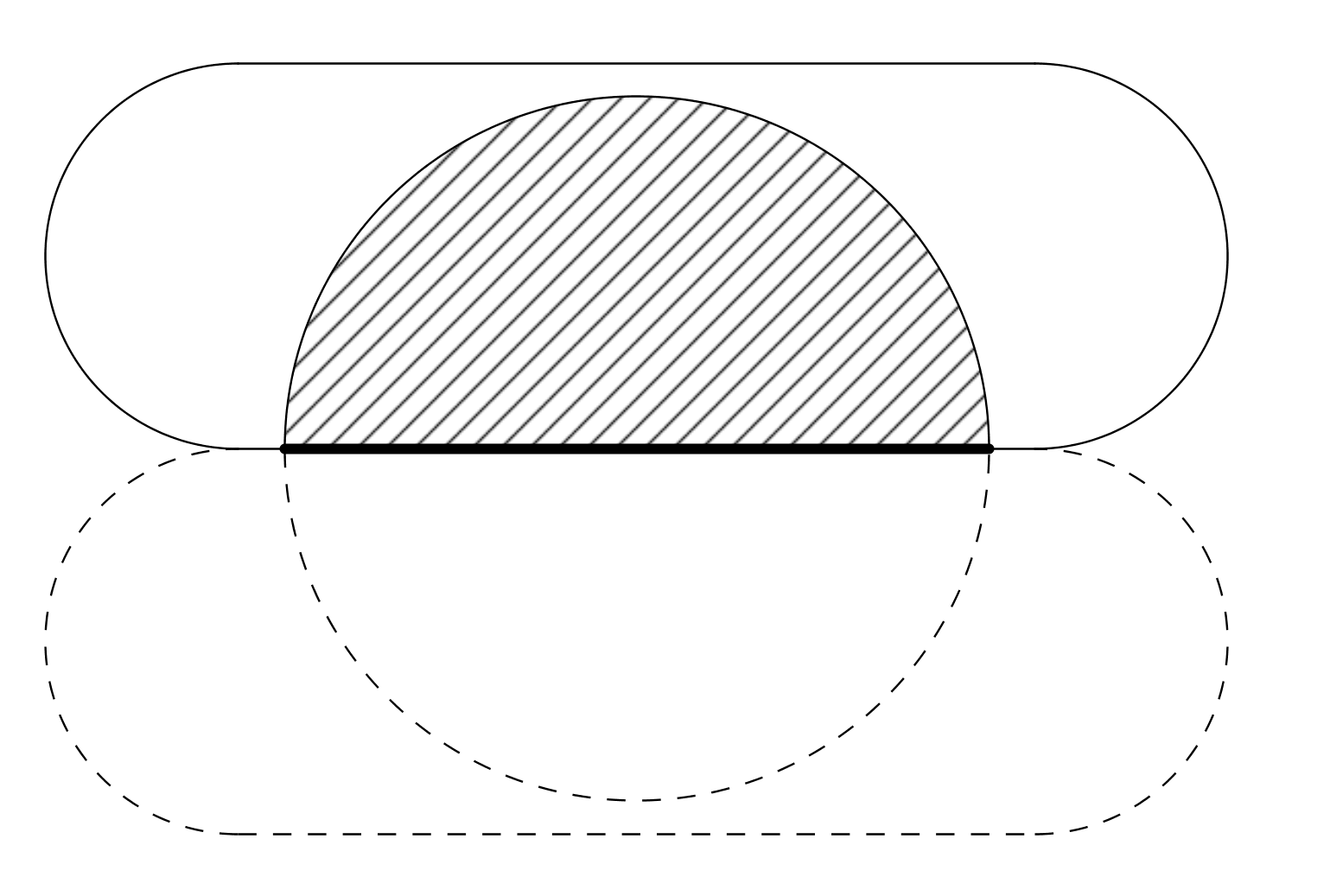}
\caption{The half disc is first embedded in a stadium-shaped domain (solid line) so that the reflecting part (thick line) lies on the boundary. It is tehcnically convenient to embed a manifold with corners to a manifold with smooth boundary so that the reflecting part still lays on the boundary. Then one can glue another copy of the larger manifold (dashed line) along a part of the boundary. The original boundary with corners is thus reflected to a manifold with smooth boundary -- the whole disc in this case.~\cite{I:refl}}%
\label{fig:refl-constr}%
\end{figure}

If~$\partial\tilde\Omega$ is to be smooth, we need to have a corner on~$\partial\Omega$ where~$E$ and~$R$ meet.
The need of corners is easily seen when one reflects a half disc (which has corners) into a whole disc (which has smooth boundary).
This example is depicted in figure~\ref{fig:refl-constr}.

Smoothness of the metric on the reflected manifold~$\tilde\Omega$ requires that~$R$ is very flat.
If the metric~$g$ on~$\Omega$ is~$C^\infty$, then the reflected metric~$\tilde g$ is~$C^1$ if and only if the second fundamental form vanishes identically on~$R$.
By symmetry $\tilde g\in C^1$ if and only if $\tilde g\in C^2$.
Higher order flatness is required on~$R$ if higher order differentiability of~$\tilde g$ is required.
If~$R$ is strictly concave or convex, solutions to the geodesic equation on~$\tilde\Omega$ branch or do not exist.
For more details on the effect of regularity of~$R$, see~\cite[Lemma~13]{I:refl}.

This reflection method is well known in the study of billiards~\cite{book-billiards}, but has only recently been used for the broken ray transform.
Hubenthal~\cite{H:square} used reflections to study the broken ray transform in the square using microlocal techniques.
He has recently~\cite{H:brt-flat-refl} used a similar methods for general planar domains where~$R$ is flat but may contain many components.
In this setting all broken rays in some neighborhood of a fixed one have the same reflection pattern (they hit the same components of~$E$ in the same order) and the domain can be unfolded so that each of these broken rays becomes a line -- note that for different patterns we obtain a different reflected domain~$\tilde\Omega$.

The main theorem~\cite[Theorem~16]{I:refl} states that if the X-ray transform is injective on~$\tilde\Omega$, then the broken ray transform is injective on~$\Omega$.
However, for the X-ray transform to make sense and to be injective on~$\tilde\Omega$, it is useful to have~$R$ flat.

Curved~$R$ is at present untractable, but corners are not.
Treating corners on Riemannian manifolds is technically difficult, but we can give Euclidean examples.
Proposition~\cite[Proposition~6]{I:refl} states that the broken ray transform is injective in a domain~$\Omega$ contained in a planar cone with any opening angle where~$R$ is contained in the boundary of the cone.
One example of such a domain is a square where~$E$ consists of two adjacent sides.

\begin{figure}%
\includegraphics{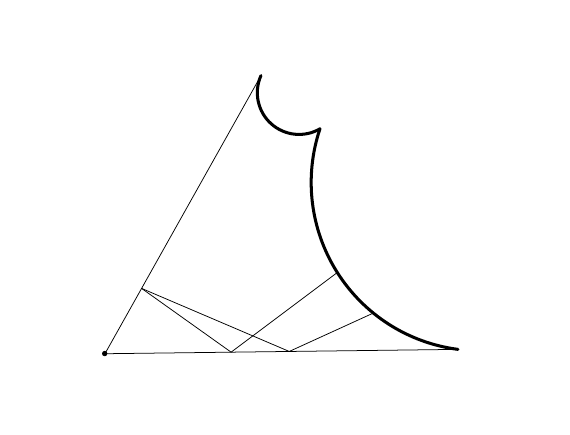}
\caption{A cone $\Omega$ with opening angle $\pi/3$ and a broken ray.~\cite{I:refl}}%
\label{fig:refl1}%
\end{figure}

\begin{figure}%
\includegraphics{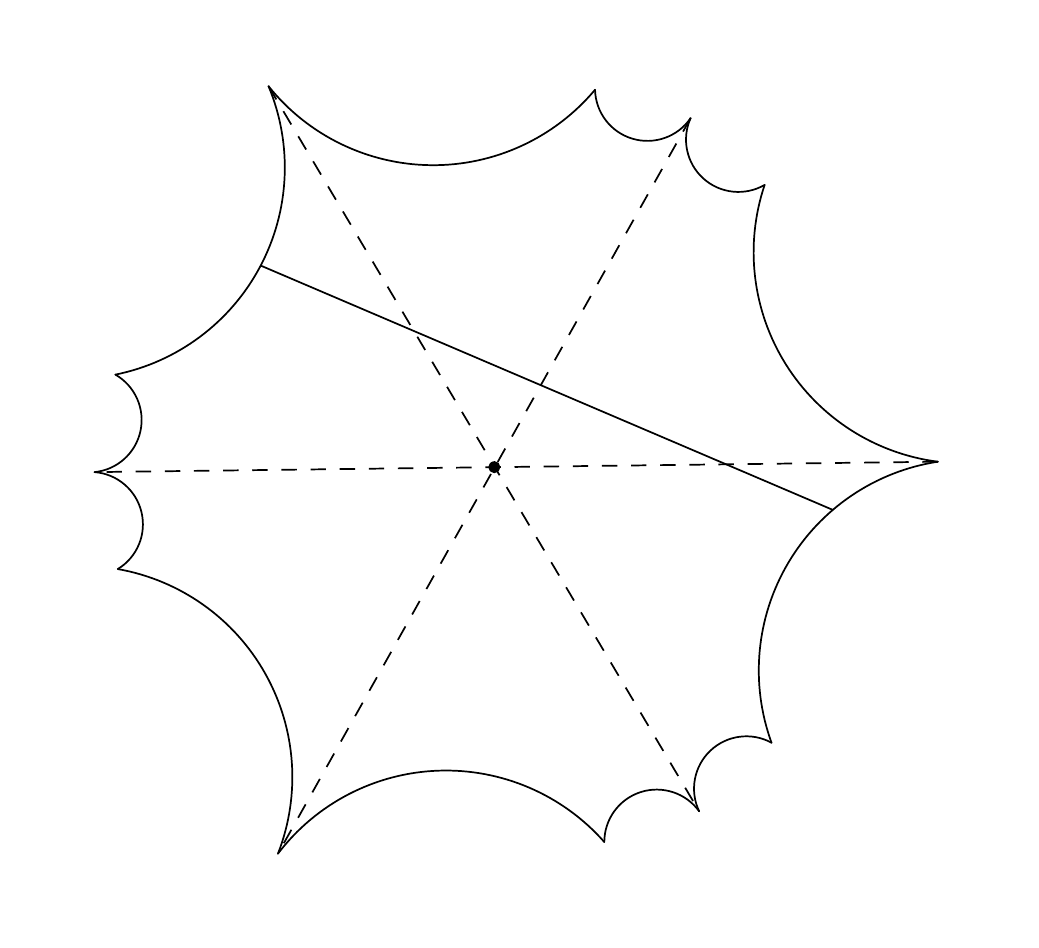}
\caption{A domain $\tilde\Omega$ obtained from the domain $\Omega$ in figure~\ref{fig:refl1} by gluing six copies together. The broken ray unfolds into a straight line.~\cite{I:refl}}%
\label{fig:refl2}%
\end{figure}

If the opening angle is $\pi/m$ for some $m\in\N$, then~$2m$ copies of~$\Omega$ can be neatly glued together (along matching sides) to form the domain~$\tilde\Omega$.
The X-ray transform is injective in planar domains, and injectivity of the broken ray transform in~$\Omega$ follows.
For general angles the result reduces to Helgason's support theorem~\cite{book-helgason} instead of injectivity of the X-ray transform.
We glue together so many copies of~$\Omega$ that~$\tilde\Omega$ has opening angle at least~$\pi$.
Then the integral of the reflected function~$\tilde f$ is known over all lines in the plane that evade a convex cone, and by the support theorem~$\tilde f$ is determined outside the cone.
These two constructions are illustrated in~\cite[Figures~2--5]{I:refl}.
The simpler case is illustrated here in figures~\ref{fig:refl1} and~\ref{fig:refl2}.

The reflection argument can also be used for the periodic broken ray transform.
Examples of this are given in~\cite[Propositions~30 and~31]{I:refl}, and here we mention the example of the square.
Just like in the above construction for cones, we take four copies of the square and glue them together.
Then we identify the opposite sides of the square -- we have thus constructed a flat torus.
Now any periodic geodesic on the torus corresponds to a periodic broken ray in the square and we arrive at a new question: do the integrals of a function on the torus over all periodic geodesics determine the function?
The answer to this question is affirmative~\cite{AR:radon-torus,I:torus}, whence the periodic broken ray transform on the square is injective.

\subsection{Boundary determination}
\label{sec:bdy-det}

%

Let $\sigma:[0,L]\to\partial\Omega$ be a geodesic on~$\partial\Omega$ with endpoints in~$\overline{\sisus E}$.
We want to approximate~$\sigma$ with a sequence~$(\gamma_n)$ of broken rays in~$\Omega$.
If $(\gamma_n(0),\dot\gamma_n(0))\to(\sigma(0),\dot\sigma(0))$, the broken rays are closer and closer to being tangent to~$\sigma$ at the starting point.
To have the broken rays remain almost tangent to~$\sigma$ at later times,~$\partial\Omega$ has to be strictly convex along~$\sigma$ -- this means that $\tpm{\dot\sigma}>0$, where $\tpm{\cdot}$ is the second fundamental form.
We call such boundary geodesics~$\sigma$ admissible.
Under these conditions $(\gamma_n,\dot\gamma_n)\to(\sigma,\dot\sigma)$ uniformly on~$[0,L]$ (see~\cite[Lemma~7]{I:bdy-det}).
If~$\partial\Omega$ is strictly concave along~$\sigma$, the broken rays~$\gamma_n$ escape to the interior of~$\Omega$.

Assuming that $f:\bar\Omega\to\C$ is continuous, this convergence implies that we may reconstruct the integral of~$f$ over admissible boundary geodesics from its broken ray transform.
If there is an admissible boundary geodesic through each point on~$E$, it is easy to reconstruct~$f$ on~$E$.
If the X-ray transform (restricted to admissible geodesics) is injective on $\bar R=\partial\Omega\setminus\sisus E$, we may thus reconstruct~$f|_R$ from the broken ray transform of~$f$.
Thus under favourable circumstances we may reconstruct~$f|_{\partial\Omega}$ from~$\brt{}f$.

This boundary determination method becomes more likely to work when the number of admissible boundary geodesics increases.
If we assume~$\partial\Omega$ to be strictly convex, every boundary geodesic is admissible.
There are some results for injectivity of the X-ray transform on Riemannian manifolds (see~\cite[p.~7]{I:bdy-det}), but current results are insufficient in very general situations.

Somewhat surprisingly, this method allows to construct not only the values of the unknown function at the boundary, but also its normal derivatives of any order.
For $k\geq1$, knowledge of lower order derivatives and the broken ray transform allows one to reconstruct the integral of $S\partial\Omega\ni (x,v)\mapsto\tpm{v}^{-k/3}\partial_\nu^kf(x)$ over all admissible geodesics.
Thus, if the X-ray transform on~$\partial\Omega$ is injective when weighted with powers of the second fundamental form, one may recover the Taylor polynomial of a function at the boundary from its broken ray transform.
This boundary determination result is given in~\cite[Corollary~2]{I:bdy-det}.

Like the approach presented in section~\ref{sec:explicit}, this one is based on the limit of infinitely many reflections.
Also, similarly to the approach in section~\ref{sec:refl}, this one reduces the problem to the X-ray transform -- this time with weight.

\subsection{PDE approach}
\label{sec:pde}

This approach to the broken ray transform was introduced by Mukhometov~\cite{M1,M2,M3,M4,M5,M6}, and we describe the main idea.
Let $\Omega\subset\R^2$ be an annular domain, that is, a strictly convex bounded domain from which a strictly convex obstacle has been removed.
Let the outer boundary be the set of tomography~$E$ and the inner boundary reflective (the set~$R$).
Suppose $f\in C^2(\bar\Omega)$ has vanishing broken ray transform.
Our goal is now to show that $f=0$.
The results in~\cite{IS:brt-pde-1obst} hold for Riemannian surfaces, but for the sake of simplicity we assume here~$\Omega$ to be Euclidean.

Let~$S\Omega$ denote the unit sphere bundle of~$\Omega$, which in the Euclidean case is simply $S\Omega=\Omega\times S^1$.
A point in~$S\Omega$ defines a point in~$\Omega$ and a direction at that point.
For $(x,v)\in S\Omega$ we have a unique broken ray $\gamma_{x,v}:[0,\tau_{x,v}]\to\bar\Omega$ satisfying $\gamma_{x,v}(0)=x$ and $\dot\gamma_{x,v}(0)=v$.
Here~$\tau_{x,v}$ is the time it takes for the broken ray to reach the set~$E$.

For a sufficiently smooth function $u:S\Omega\to\R$ we define the derivatives (vector fields)~$X$ and~$V$ by letting $Xu(x,v)=v\cdot\nabla_xu(x,v)$ and $Vu(x,v_\theta)=\partial_\theta u(x,v_\theta)$, where we have written $v_\theta=(\cos\theta,\sin\theta)$.

Now we define a function~$u^f$ on~$S\Omega$ by letting
\begin{equation}
u^f(x,v)
=
\int_0^{\tau_{x,v}}f(\gamma_{x,v}(t))\der t.
\end{equation}
By the fundamental theorem of calculus this function satisfies $Xu^f(x,v)=-f(x)$.
But since~$f(x)$ does not depend on the direction~$v$, we get $VXu^f(x,v)=0$ for all $(x,v)\in S\Omega$.

We would like the equation $VXu=0$ to have a unique solution, but this obviously requires that we fix some boundary conditions.
The boundary of $S\Omega$ is $E\times S^1\cup R\times S^1$.
The function~$f$ has vanishing broken ray transform, so $u^f(x,v)=0$ for all $x\in E$.
The boundary condition on~$R$ is more complicated.

At every $x\in R$ we have the outer unit normal~$\nu(x)$.
We define the map $\rho_x:S^1\to S^1$ by $\rho_x v=v-2\ip{v}{\nu(x)}\nu(x)$; this map reflects~$v$ accross the tangent of~$R$ at~$x$.
Since broken rays reflect so that the incoming and outgoing directions~$v_i$ and~$v_o$ are related via $v_o=\rho_x v_i$ (for a reflection at~$x$), we have that $u^f(x,v)=u^f(x,\rho_x v)$ for all $x\in R$.
This is the boundary condition on~$R$.

It thus remains to show that $u=0$ is the only solution to
\begin{equation}
\label{eq:PDE}
\begin{cases}
VXu=0 & \text{in }S\Omega\\
u=0 & \text{on }E\times S^1\\
u(x,v)=u(x,\rho_xv) & \text{on }R\times S^1.
\end{cases}
\end{equation}
We cannot assume that $u\in C^2(S\Omega)$ even if $f\in C^2(\bar\Omega)$, since broken rays do not depend smoothly on their initial point and direction close to tangential reflections.
We will return to this regularity issue shortly.

As it turns out, for $u\in C^2(S\bar\Omega)$ we have~\cite[Lemma~8]{IS:brt-pde-1obst}
\begin{equation}
\label{eq:vv1}
\begin{split}
\aabs{VXu}_{L^2(S\Omega)}^2
&=
\aabs{XVu}_{L^2(S\Omega)}^2
+
\aabs{Xu}_{L^2(S\Omega)}^2
\\&\quad
-
\iip{KVu}{Vu}_{L^2(SM)}
+
\iip{\nabla_Tu}{Vu}_{L^2(\partial SM)}
,
\end{split}
\end{equation}
where $\nabla_T=\ip{\nu^\perp}{\nabla}$ is the tangential derivative and~$K$ is the Gaussian curvature of~$\Omega$.
Of course $K=0$ for an Euclidean domain~$\Omega$, but the same identity holds true on Riemannian surfaces where~$K$ need not vanish.

One can further manipulate the boundary term of equation~\eqref{eq:vv1}.
If~$\kappa$ denotes the signed curvature of~$\partial M$ and~$u_e$ and~$u_o$ are the even and odd part of~$u$ on~$\partial M$ with respect to the reflection~$\rho$, we have~\cite[Lemma~9]{IS:brt-pde-1obst}
\begin{equation}
\label{eq:vv2}
\begin{split}
\iip{\nabla_Tu}{Vu}_{L^2(\partial SM)}
&=
\iip{\nabla_Tu_e}{Vu_o}_{L^2(\partial SM)}
+
\iip{\nabla_Tu_o}{Vu_e}_{L^2(\partial SM)}
\\&\quad-
\iip{\kappa Vu}{Vu}_{L^2(\partial SM)}
.
\end{split}
\end{equation}
Combining the identities~\eqref{eq:vv1} and~\eqref{eq:vv2} with the boundary conditions $u=0$ on $E\times S^1$ and $u_o=0$ on $R\times S^1$, we find
\begin{equation}
\label{eq:pestov}
\begin{split}
\aabs{VXu}_{L^2(S\Omega)}^2
&=
\aabs{XVu}_{L^2(S\Omega)}^2
+
\aabs{Xu}_{L^2(S\Omega)}^2
\\&\quad
-
\iip{KVu}{Vu}_{L^2(SM)}
-
\iip{\kappa Vu}{Vu}_{L^2(\partial SM)}.
\end{split}
\end{equation}
Energy estimates of this type are known as Pestov identities.

We need~$u$ to have enough regularity so that the calculations leading to the identity~\eqref{eq:pestov} can be justified.
Proving this regularity is rather involved, and we omit it here.
The outcome of the regularity analysis is that if $f\in C^2(\bar\Omega)$ has vanishing broken ray transform, then $u^f$ has enough regularity for the identity~\eqref{eq:pestov} to hold.

The obstacle is convex, so $\kappa\leq0$.
In Euclidean domains $K=0$, and more generally on nonpositively curved Riemannian surfaces $K\leq0$.
These conditions give a positive sign to the last two terms in~\eqref{eq:pestov}.
As observed above, $VXu^f=0$ and $Xu^f=-f$, so~\eqref{eq:pestov} yields
\begin{equation}
0
\geq
\aabs{f}_{L^2(\Omega)}^2.
\end{equation}
This implies that $f=0$, which we set out to prove.

If we could prove the Pestov identity~\eqref{eq:pestov} for~$u^f$ without assuming that the broken ray transform of $f\in C^2(\bar\Omega)$ vanishes, we would obtain
\begin{equation}
C\aabs{Vu^f}^2_{L^2(E\times S^1)}
\geq
\aabs{f}^2_{L^2(\Omega)},
\end{equation}
where the constant~$C$ depends on the maximum of the curvature~$\kappa$ at~$E$.
Since with natural identifications $u^f|_{E\times S^1}=\brt f$, this would be a stability estimate for the broken ray transform.

In the Euclidean case this result actually follows easily from Helgason's support theorem~\cite{book-helgason}.
The support theorem states that if a continuous function in the bounded domain~$\Omega$ integrates to zero over all lines that avoid the convex obstacle, then it has to vanish outside the obstacle.
Similar support theorems are also available on manifolds~\cite{K:spt-thm,UV:local-x-ray}, but they cannot be used on two dimensional manifolds with nonanalytic metric.

Thus the result in~\cite{IS:brt-pde-1obst} only gives new information on nonanalytic surfaces.
The approach used is suitable for tackling the same problem with several convex obstacles, whereas support theorems can only give information of the exterior of the convex hull of the union of the obstacles.

\begin{figure}%
\includegraphics[scale=.7]{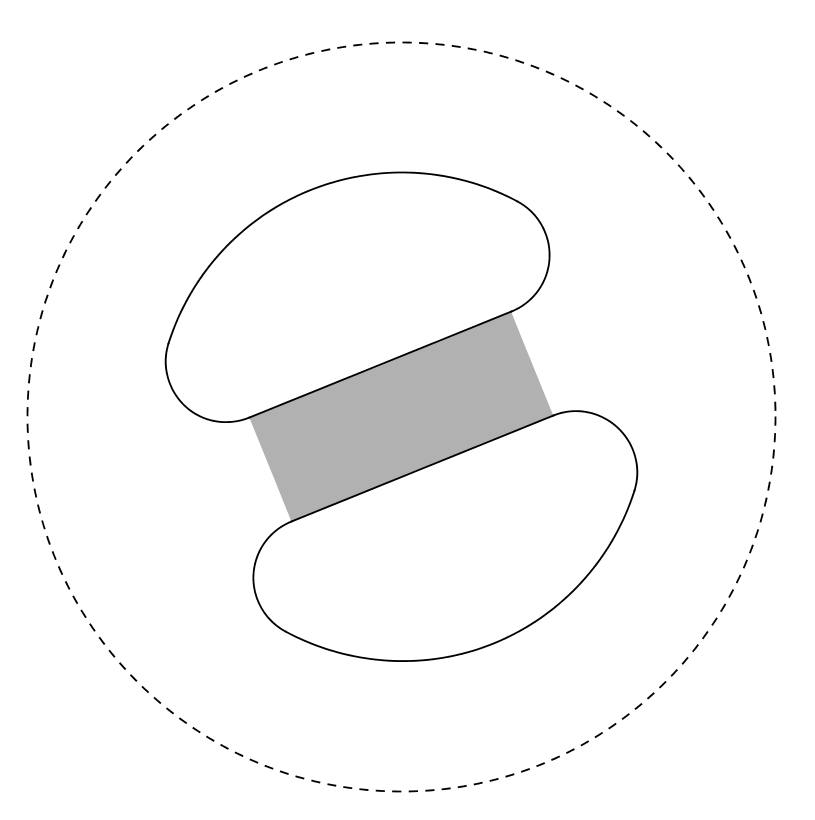}
\caption{A planar domain with two obstacles where the broken ray transform fails to be injective. One can construct nonzero functions supported in the gray area such that the broken ray transform vanishes.~\cite{IS:brt-pde-1obst}}%
\label{fig:2obst}%
\end{figure}

The broken ray transform was shown to be injective in the case of several reflecting obstacles in the Euclidean plane by Eskin~\cite{eskin}.
That result, however, relies on stringent geometrical assumptions on the obstacles, which for example prevent them from being smooth.
Showing injectivity of the broken ray transform in a planar Euclidean domain with two smooth, convex obstacles remains an open problem.
Actually, if one only assumes the two obstacles to be convex, there is an example (see figure~\ref{fig:2obst}) of noninjective broken ray transform, so at least one of the obstacles should be strictly convex.

Pestov identities have been widely used in the study of the X-ray transform.
Such energy identities have proven to be useful when there is not enough structure or symmetry for explicit calculations.
For a review of the use of Pestov identities, see~\cite[Section~4]{PSU:tensor-survey}.
Mukhometov has also used Pestov identities to study the broken ray transform~\cite{M1, M2, M3, M4, M5, M6}.

\section{Examples and counterexamples}
\label{sec:(ctr)ex}

To make the results more concrete, we list examples of domains~$\Omega$ (with the corresponding set of tomography~$E$ or its complement $R=\partial\Omega\setminus E$) where the broken ray transform has been proven to be injective:
\begin{itemize}
\item Euclidean ball with arbitrarily small open~$E$ for functions that are analytic in the angular variable in a suitable sense~\cite{I:disk},
\item any Euclidean cone for piecewise continuous functions where~$R$ is contained in the surfaces of the cone~\cite{I:refl},
\item quarter of~$S^n$, $n\geq3$, where~$E$ is half of the boundary for compactly supported~$L^2$ functions~\cite{I:refl},
\item two dimensional hemisphere where~$E$ is slightly larger than half of the boundary for compactly supported smooth functions~\cite{I:refl}, and
\item nonpositively curved compact Riemannian surface with strictly convex boundary~$E$ with a convex obstacle removed (boundary of the obstacle is~$R$) for~$C^2$ functions~\cite{IS:brt-pde-1obst}.
\end{itemize}

Suppose we know the broken ray transform of a function in the Euclidean ball (of dimension three or higher) which is smooth in a neighborhood of the boundary.
Assuming that~$E$ is open and contains a hemisphere, we may reconstruct the Taylor polynomial of the function at all boundary points from its broken ray transform.~\cite[after Remark~5]{I:bdy-det}

There are also examples for the periodic broken ray transform, where all of the boundary is reflective.
These examples are an octant of~$S^2$ for compactly supported smooth functions~\cite{I:refl} and the cube $[0,1]^n$, $n\geq2$, for a sum of an~$L^1$ function and a compactly supported distribution~\cite{I:torus}.

In addition to examples, we also have counterexamples.
If a manifold with boundary contains a (generalized) reflecting tubular part, then the broken ray transform (periodic or regular) is not injective~\cite{I:refl}.
For a precise statement see~\cite[Proposition~29]{I:refl}; we only remark here that this includes the counterexample of figure~\ref{fig:2obst}.
There is also a counterexample for the periodic broken ray transform; namely, the periodic broken ray transform is not injective on compactly supported smooth functions in the Euclidean unit disc~\cite{I:refl}.

Many of these examples and counterexamples are described in more detail in~\cite[Sections~6--7]{I:refl}.

\section{Applications}
\label{sec:app}

The X-ray transform has important applications in medical imaging.
Inferring the three dimensional structure of an object from X-ray images taken from all directions relies on our ability to invert the X-ray transform.
This mathematical model underlies CT, PET and SPECT imaging.
For details of these applications, see~\cite{book-natterer}.

The X-ray transform also arises as a linearization of other inverse problems.
For example, if one linearizes the boundary rigidity problem (see section~\ref{sec:ip}) with respect to a conformal variation of the metric, one ends up with the X-ray transform.
If one has a suitable form of the inverse function theorem (the spaces in question are infinite dimensional, so this is not obvious), solving the linearized problem shows that the original problem can be solved locally.
For an example of a boundary rigidity problem solved via linearization, see~\cite{SUV:partial-bdy-rig}.

Applications of the broken ray transform are described in~\cite[Section~6]{IS:brt-pde-1obst}, and we give a brief summary of that discussion here.

Imaging with broken rays is relevant for seismology.
To give a hint of this vast area of research, we refer to~\cite{G:refl-tomo,DGGS:refl-tomo,DL:refl-tomo}.

Eskin~\cite{eskin} reduced an inverse boundary value problem for the electromagnetic Schr\"odinger operator to the injectivity of the broken ray transform.
Similarly, Kenig and Salo~\cite{KS:calderon} reduced a partial data problem for Calder\'on's problem in a tubular domain to the injectivity of the broken ray transform in the transversal domain.
These results are based on constructing solutions to the corresponding PDEs which concentrate near broken rays.
Therefore new results for the broken ray transform imply results for these problems; an example is given in~\cite[Theorem~3]{I:refl}.


Like the X-ray transform, the broken ray transform also arises as a linearization.
Let $(M,g)$ be a Riemannian manifold with boundary and let its boundary be partitioned in two sets~$E$ and~$R$ as above, and suppose~$E$ is open.
Let $x\in E$ be a boundary point and $v\in S_xM$ a unit vector at it.
Assuming that~$v$ points into the interior of~$M$, there is a unique broken ray~$\gamma_{x,v}$ with initial point and velocity $(x,v)$.
Suppose $\gamma_{x,v}$ is has no tangential reflections.
Let then $f\in C^\infty(M)$ be any function and define new metrics $g_s=(1+sf)g$ for all $s\in(-\eps,\eps)$.
If~$\eps$ is small enough, the tensors~$g_s$ really are Riemannian metrics.
Then let $\gamma^s_{x,v}$ be the unique broken ray with respect to the metric~$g_s$ starting at the point~$(x,v)$.
By making~$\eps$ smaller if necessary we may assume that none of these broken rays has tangential reflections.

Let~$\tau_s$ denote the length of $\gamma^s_{x,v}$ with respect to the metric~$g_s$.
If the broken rays $\gamma^s_{x,v}$ have the same endpoint, then the derivative $\left.\Der{s}\tau_s\right|_{s=0}$ is the integral of the function~$2f$ over the broken ray~$\gamma_{x,v}$.
In this sense the broken ray transform is the linearization of the boundary distance function with reflections.
This result holds in more generality; see~\cite[Theorem~17]{IS:brt-pde-1obst}.

Ray transforms are also related to spectral geometry.
Spectral properties of a domain or a manifold can often be related to the lengths of periodic geodesics or broken rays (billiard trajectories) on it; see~\cite{DH:spectral-survey} for a review of the topic.
Since lengths of broken rays correspond to the broken ray transform via linearization as described above, it is natural to expect the periodic broken ray transform to have applications in spectral geometry.
As far as we know, the precise connection between the periodic broken ray transform and spectral geometry of manifolds with boundary is yet to be discovered.

\newpage

\bibliographystyle{abbrv}
\bibliography{brt}

\begin{thebibliography}{10}

\bibitem{AR:radon-torus}
A.~Abouelaz and F.~Rouvi\`ere.
\newblock {Radon Transform on the Torus}.
\newblock {\em Mediterranean Journal of Mathematics}, 8(4):463--471, 2011.

\bibitem{A:radon}
G.~Ambartsoumian.
\newblock {Inversion of the {V}-line {R}adon transform in a disc and its
  applications in imaging}.
\newblock {\em Computers \& Mathematics with Applications}, 64(3):260--265,
  2012.

\bibitem{AM:v-line-disc}
G.~Ambartsoumian and S.~Moon.
\newblock A series formula for inversion of the {V}-line {R}adon transform in a
  disc.
\newblock {\em Computers \& Mathematics with Applications}, 66(9):1567--1572,
  2013.

\bibitem{cormack}
A.~M. Cormack.
\newblock {Representation of a Function by Its Line Integrals, with Some
  Radiological Applications}.
\newblock {\em Journal of Applied Physics}, 34(9):2722--2727, 1963.

\bibitem{DH:spectral-survey}
K.~{Datchev} and H.~{Hezari}.
\newblock {Inverse problems in spectral geometry}.
\newblock In {\em Inverse problems and applications: Inside Out II}, volume~60
  of {\em Math. Sci. Res. Inst. Publ.}, pages 455--486. Cambridge Univ. Press,
  Cambridge, 2012.

\bibitem{DGGS:refl-tomo}
F.~Delbos, J.~C. Gilbert, R.~Glowinski, and D.~Sinoquet.
\newblock Constrained optimization in seismic reflection tomography: a
  gauss-newton augmented lagrangian approach.
\newblock {\em Geophysical Journal International}, 164(3):670--684, 2006.

\bibitem{DL:refl-tomo}
F.~Delprat-Jannaud and P.~Lailly.
\newblock Ill-posed and well-posed formulations of the reflection travel time
  tomography problem.
\newblock {\em Journal of Geophysical Research: Solid Earth},
  98(B4):6589--6605, 1993.

\bibitem{eskin}
G.~Eskin.
\newblock {Inverse boundary value problems in domains with several obstacles}.
\newblock {\em Inverse Problems}, 20(5):1497--1516, 2004.

\bibitem{G:refl-tomo}
S.~V. Goldin.
\newblock {Ray Reflection Tomography: Review and Comments}.
\newblock In {\em 3rd International Congress of the Brazilian Geophysical
  Society}, 1993.

\bibitem{book-helgason}
S.~Helgason.
\newblock {\em {The Radon Transform}}.
\newblock Birkh\"auser, 2. edition, 1999.

\bibitem{H:brt-flat-refl}
M.~{Hubenthal}.
\newblock {The Broken Ray Transform in $n$ Dimensions}.
\newblock Oct. 2013.
\newblock \href{http://arxiv.org/abs/1310.7156}{arXiv:1310.7156}.

\bibitem{H:square}
M.~{Hubenthal}.
\newblock {The Broken Ray Transform On The Square}.
\newblock Feb. 2013.
\newblock \href{http://arxiv.org/abs/1302.6193}{arXiv:1302.6193}.

\bibitem{KK:brt-inversion-range}
A.~Katsevich and R.~Krylov.
\newblock Broken ray transform: inversion and a range condition.
\newblock {\em Inverse Problems}, 29(7):075008, 2013.

\bibitem{KS:calderon}
C.~E. {Kenig} and M.~{Salo}.
\newblock {The Calder\'on problem with partial data on manifolds and
  applications}.
\newblock {\em Analysis \& PDE}, 2012.
\newblock To appear, \href{http://arxiv.org/abs/1211.1054}{arXiv:1211.1054}.

\bibitem{K:spt-thm}
V.~P. Krishnan.
\newblock {A Support Theorem for the Geodesic Ray Transform on Functions}.
\newblock {\em Journal of Fourier Analysis and Applications}, 15:515--520,
  2009.

\bibitem{KLU:broken-geodesic-flow}
Y.~Kurylev, M.~Lassas, and G.~Uhlmann.
\newblock Rigidity of broken geodesic flow and inverse problems.
\newblock {\em Amer. J. Math.}, 132(2):529--562, 2010.

\bibitem{MNTZ:radon}
M.~Morvidone, M.~K. Nguyen, T.~T. Truong, and H.~Zaidi.
\newblock {On the V-Line Radon Transform and Its Imaging Applications}.
\newblock {\em International Journal of Biomedical Imaging}, 2010(208179),
  2010.

\bibitem{M:bdy-rig-surface-eng}
R.~G. Mukhometov.
\newblock The problem of recovery of two-dimensional riemannian metric and
  integral geometry.
\newblock {\em Soviet Math. Dokl.}, 18(1):27--31, 1977.

\bibitem{M:bdy-rig-surface-rus}
R.~G. Mukhometov.
\newblock The problem of recovery of two-dimensional riemannian metric and
  integral geometry.
\newblock {\em Dokl. Akad. Nauk SSSR}, 232(1):32--35, 1977.

\bibitem{M1}
R.~G. Mukhometov.
\newblock A problem of integral geometry on a plane with respect to a family of
  wave fronts ({R}ussian).
\newblock In {\em Methods for investigating ill-posed problems of mathematical
  physics}, pages 56--68, 103. Akad. Nauk SSSR Sibirsk. Otdel., Vychisl.
  Tsentr, Novosibirsk, 1983.

\bibitem{M2}
R.~G. Mukhometov.
\newblock A problem of integral geometry in a class of compactly supported
  functions ({R}ussian).
\newblock In {\em Linear and nonlinear problems in computer tomography}, pages
  124--131. Akad. Nauk SSSR Sibirsk. Otdel., Vychisl. Tsentr, Novosibirsk,
  1985.

\bibitem{M3}
R.~G. Mukhometov.
\newblock {\em Estimates for the stability of solutions of problems of integral
  geometry in the case of nonregular families of curves ({R}ussian)}.
\newblock Preprint, volume 681, Akad. Nauk SSSR Sibirsk. Otdel., Vychisl.
  Tsentr, Novosibirsk, 1986.

\bibitem{M4}
R.~G. Mukhometov.
\newblock Problems of integral geometry in a domain with a reflecting part of
  the boundary {{(Russian, translation in Soviet Math. Dokl. 36 (1988), no. 2,
  260--264)}}.
\newblock {\em Dokl. Akad. Nauk SSSR}, 296(2):279--283, 1987.

\bibitem{M5}
R.~G. Mukhometov.
\newblock A problem of integral geometry for a family of rays with multiple
  reflections.
\newblock In {\em Mathematical methods in tomography ({O}berwolfach, 1990)},
  volume 1497 of {\em Lecture Notes in Math.}, pages 46--52. Springer, Berlin,
  1991.

\bibitem{M6}
R.~G. Mukhometov.
\newblock On problems of integral geometry in the non-convex domains.
\newblock In {\em Tomography, impedance imaging, and integral geometry ({S}outh
  {H}adley, {MA}, 1993)}, volume~30 of {\em Lectures in Appl. Math.}, pages
  157--176. Amer. Math. Soc., Providence, RI, 1994.

\bibitem{book-natterer}
F.~Natterer.
\newblock {\em {The Mathematics of Computerized Tomography}}.
\newblock Society for Industrial and Applied Mathematics, 2001.

\bibitem{PSU:tensor-survey}
G.~P. {Paternain}, M.~{Salo}, and G.~{Uhlmann}.
\newblock {Tensor tomography: progress and challenges}.
\newblock {\em Chinese Ann. Math. Ser. B}, 2013.
\newblock To appear, \href{http://arxiv.org/abs/1303.6114}{arXiv:1303.6114}.

\bibitem{radon}
J.~Radon.
\newblock {\"Uber die Bestimmung von Funktionen durch ihre Integralwerte
  l\"angs gewisser Mannigfaltigkeiten}.
\newblock {\em {Berichte der Sachsischen Akadamie der Wissenschaft}},
  69:262--277, 1917.

\bibitem{S:tensor-book}
V.~A. Sharafutdinov.
\newblock {\em Integral geometry of tensor fields}.
\newblock Walter de Gruyter, 1994.

\bibitem{SUV:partial-bdy-rig}
P.~{Stefanov}, G.~{Uhlmann}, and A.~{Vasy}.
\newblock {Boundary rigidity with partial data}.
\newblock June 2013.
\newblock \href{http://arxiv.org/abs/1306.2995}{arXiv:1306.2995}.

\bibitem{book-billiards}
S.~Tabachnikov.
\newblock {\em {Geometry and Billiards}}, volume~30 of {\em Student
  Mathematical Library}.
\newblock American Mathematical Society, Providence, RI, 2005.

\bibitem{TN:radon}
T.~T. Truong and M.~K. Nguyen.
\newblock On new {$\mathfrak V$}-line radon transforms in $\mathbb{R}^{2}$ and
  their inversion.
\newblock {\em Journal of Physics A: Mathematical and Theoretical},
  44(7):075206, 2011.

\bibitem{UV:local-x-ray}
G.~{Uhlmann} and A.~{Vasy}.
\newblock {The inverse problem for the local geodesic ray transform}.
\newblock Oct. 2012.
\newblock \href{http://arxiv.org/abs/1210.2084}{arXiv:1210.2084}.

\bibitem{ZSM:star-transform}
F.~{Zhao}, J.~C. {Schotland}, and V.~A. {Markel}.
\newblock {Fourier-space inversion of the star transform}.
\newblock Jan. 2014.
\newblock \href{http://arxiv.org/abs/1401.7655}{arXiv:1401.7655}.

\end{thebibliography}


\begin{thebibliography}{999}
\bibitem[I]{I:disk}
{\bf Broken ray tomography in the disc} \\
J.~Ilmavirta, \\
Inverse Problems  29(3):035008 (2013).
\\~
\bibitem[II]{I:refl}
{\bf A reflection approach to the broken ray transform} \\
J.~Ilmavirta, \\
to appear in Mathematica Scandinavica.
\\~
\bibitem[III]{I:bdy-det}
{\bf Boundary reconstruction for the broken ray transform} \\
J.~Ilmavirta, \\
Annales Academiae Scientiarum Fennicae Mathematica 39(2):485--502 (2014).
\\~
\bibitem[IV]{I:torus}
{\bf On Radon transforms on tori} \\
J.~Ilmavirta, \\
preprint.
\\~
\bibitem[V]{IS:brt-pde-1obst}
{\bf Broken ray transform on a Riemann surface with a convex obstacle} \\
J.~Ilmavirta, M.~Salo, \\
preprint.
\\~
\end{thebibliography}

\newpage


\begin{appendices}
\begin{sloppypar}
\section{A beginner's introduction to inverse problems}
\label{sec:eng}

One constantly encounters situations where something has to be measured: speed of a car, volume of flour, power of a lamp, or something else.
It is of utmost importance that suitable tools are available for measurement so that one may successfully measure -- and live.

Measurement requires two kinds of tools, technical and mental.
If we want to measure the speed of a car, we drive it on a road, measure the spent time by a clock, and count the number of times a wheel has gone around to figure out the distance.
We also need to know the circumference of the wheel, which we can easily measure with a tape measure.
Thus, using three technical tools (clock, counter and tape measure), we have found three quantities, but none of them is the speed.

Next we need mental tools, that is, methods for reasoning and computing.
Multiplying the circumference by the number of rotations gives the total length.
Dividing this by the time gives the average speed.
We have now finally measured the speed, and both technical and mental tools have been of use.

Sometimes a measurement cannot be done directly.
If a physician wants to measure the density of a patient's thighbone, he should not remove the bone from the thigh and measure its volume and mass to obtain the density as their quotient.
Care for the patient's health poses a limit for the method, and the measurement is to be done indirectly.
One has to use some other way to figure out the density of the bone.
Because this reasoning is not always easy, the importance of mental tools is great in comparison to technical ones in indirect measurements.

\subsection{Direct and inverse problems}

The goal of inverse problems research is to produce tools of both kinds for indirect measurements.
My own research is focused on the mathematical side of inverse problems, so for me inverse problems are somewhat synonymous with the mathematics of indirect measurements.
Mental tools become central, and multiplication and division no longer suffice for mathematical tools.

The nature of a direct problems is the following: when the properties of an object are known, one has to determine its behaviour.
If for example the length, thickness, density and tension of a guitar string are known, one can figure out the sound it makes.
In this case the sound is composed of a fundamental frequency and its multiples.
From the mentioned properties one can determine all the possible frequencies, but not their strengths.

An inverse problem asks the same question in a different direction: when the behaviour is known, one has to determine the properties.
In the guitar example we would hear a guitar string, and from the sound we should figure out the properties of the string (length, thickness, density and tension).
But this cannot be done, since if for example the diameter is halved and length doubled, the sound stays the same.
The inverse problem does not have a clear solution, but something is known about the properties of the string.

If three of the four properties of the string are known and the frequencies have been measured, the fourth one can be figured out.
If we can measure the length and thickness and we know the material (we can check a table book for its density), we can work out the tension.
When formulated in this way, the inverse problem does indeed have a unique solution.

What happens if the guitar is replaced with a drum?
If the properties of the membrane of the drum are known, one can calculate the sound it makes (direct problem).
But if we hear a drum and we known all about the membrane but its shape, can we deduce the shape?
No one can give full answer to this question, but something is known.

If the sound sounds like it could come from a perfectly round drum, then it certainly comes from a round drum.
In other words, one can always recognize a round drum by its sound.
But there are examples of different drums with corners which make an identical sound, so there are also drums that cannot be recognized.

What drums can then be recognized by their sound?
This problem is under intense study, and the problems are in the mental tools.
This problem can be formulated mathematically, but solving it is very difficult.
Eagerness to this research does not stem from a widespread drum hobby, but the same mathematical model has also other applications.

The direct and inverse problems described above represent a field called integral geometry.
As its name suggests, it studies the relation between spectrum (frequencies of a sound) and geometry (length or shape).
There are also other kinds of inverse problems, but these make a good example.

An example of an indirect measurement resembling the instrument examples is known to many.
Water melon farmers can tell whether a melon is ripe by thumping it and listening, and similar reasoning has been used for ages.

Inverse problems are laborious to study.
As the examples demonstrate, a problem does not always have a unique solution (the same behaviour can be caused by many different properties) and a problem is more difficult than the corresponding direct problem.
It is very interesting -- and useful for applications -- only to find out, which inverse problems have unique solutions, or in other words, which methods of indirect measurement could work even in principle.

\subsection{Weighing a rope}

Let us look for a better understanding of indirect measurements by means of an example which is related to the previous ones.
We set out to weigh a rope which is fixed to a wall from both ends like a slack clothesline.
If we can detach the rope and put it on a scale, the problem is easy, but it becomes more difficult if we disallow detaching it.
Such prohibition from breaking the object under study is more natural in medical applications, because a broken patient cannot be rebuilt, but let us not be distracted by this.

One possibility is to use the rope as a guitar string and use earlier ideas, and so we shall do.
We place a bar of some kind next to one fastening point and let the rope lay over it.
(If we cannot find a suitable support for the bar, we ask a friend to hold it still.)
In the short part of the rope between this support and the wall we place a weight whose mass is known and much greater than that of the rope.

Now the longer part of the rope is tense and can be played as a guitar string.
Pulling the rope slightly down in the middle and letting go, we can make it oscillate beautifully.
We measure the time needed for ten oscillations by a clock to find out the frequency.
If all of the rope goes up and down in the same rhythm, this frequency is the fundamental frequency.

The length and thickness of the rope are easy to measure.
In addition we need to know the tension.
That we can find out, when we know the mass of the weight and measure the angle of the rope on which it hangs.

Now we are in the situation described earlier: we know the length, thickness and tension of the rope and its fundamental frequency.
From these we can calculate is density, and multiplying it with the volume gives the mass.
This is of course only the mass of the tense part, but the mass of the entire rope is easily figured out since we can see what portion of the rope is tense.
We have thus finally found the mass of the rope.

\subsection{X-ray tomography}

The solution of one inverse problem has made the life of many people easier and deserves a mention.
But let us begin with the direct problem and description of the situation.

X-rays are very similar to light rays.
X-ray radiation is light whose color the human eye cannot see -- it is in a way too blue.
Such light can be produced and measured with a suitable lamp and camera, so a properly equipped human can see with X-rays.
Unlike usual light it mostly goes through a human body, so it can be used to see inside a human.

A broken bone is easy to see in an X-ray image, but some structures are difficult to see in one image.
And even if there were multiple images, it can be difficult to see something inside a bone, for example, because something else is always in the way.

If one knows the precise structure of a human body, that is, what is in inside and where exactly it is, one can fairly easily reason what an X-ray image taken from any direction should look like.
The corresponding inverse problem asks whether one could figure out the exact three dimensional structure of the body from X-ray images from all directions.
It turns out that this is indeed possible, but it is not easy.
One needs the help of a computer because the procedure requires complicated calculations, but this is easy for a modern computer.
Due to this dependency on computers this method is often called computerized tomography (CT).

This indirect measurement method is very useful.
Every day in hospitals around the world doctors use it to see inside their patients.
It is therefore no wonder that it earned its discoverers Cormack and Hounsfield the Nobel prize in medicine in 1979.

My own research is related to inverse problems similar to X-ray tomography, but describing in more detail does not fit here.
However, inverse problems of different kinds have suprising connections, and for example the inverse problem of X-ray tomography is closely related to spectral geometry which was described above.
This is the salt of mathematical research; hardly anyone could have guessed that hearing the shape of a drum is eventually very similar to finding a fracture in a bone!

\newpage

\begin{otherlanguage}{finnish}
\section{Aloittelijan johdatus inversio-ongelmiin}
\label{sec:fin}

Jatkuvasti tulee vastaan tilanteita, joissa jotain pit‰‰ mitata: auton vauhtia, jauhojen tilavuutta, lampun tehoa ja milloin mit‰kin.
On ehdottoman t‰rke‰‰, ett‰ mittaamiseen on k‰ytˆss‰ siihen sopivia v‰lineit‰, jotta mittaaminen -- ja sen avulla el‰m‰ -- sujuu.

Mittaamisessa tarvitaan kahdenlaisia v‰lineit‰, teknisi‰ ja henkisi‰.
Jos haluamme mitata auton vauhdin, ajamme jonkin tienp‰tk‰n ja mittaamme ajamiseen kuluneen ajan kellolla ja ajomatkan laskemalla montako kierrosta auton rengas on matkan aikana pyˆr‰ht‰nyt.
Lis‰ksi tarvitsemme tiedon renkaan ymp‰rysmitasta, ja sen voimme tehd‰ mittanauhalla.
N‰in teknisi‰ apuv‰lineit‰ (kello, kierroslaskuri ja mittanauha) k‰ytt‰en olemme saaneet kolme mittaustulosta, mutta mik‰‰n n‰ist‰ ei ole auton vauhti.

Seuraavaksi tarvitaan henkisi‰ apuv‰lineit‰, nimitt‰in p‰‰ttely‰ ja laskutaitoa.
Kertomalla kierrosm‰‰r‰ renkaan ymp‰rysmitalla saadaan auton kulkema matka.
Jakamalla t‰m‰ matka kuluneella ajalla saadaan auton keskim‰‰r‰inen vauhti.
N‰in vauhti on lopulta saatu mitattua, ja siin‰ tulivat tarpeeseen niin tekniset kuin henkisetkin v‰lineet.

Joskus mittausta ei voida tehd‰ suoraan.
Jos l‰‰k‰ri haluaa mitata potilaansa reisiluun tiheyden, h‰nen ei kannata poistaa luuta potilaasta ja sen j‰lkeen mitata sen tilavuutta ja massaa, joiden osam‰‰r‰n‰ h‰n saisi tiheyden.
Potilaan terveyden vaaliminen asettaa mittaustavalle rajoitteen, ja mittaus on teht‰v‰ ep‰suorasti.
On jotenkin muuten kyett‰v‰ p‰‰ttelem‰‰n, mik‰ luun tiheys on.
Koska t‰m‰ p‰‰ttely ei ole aina helppoa, korostuu ep‰suorassa mittauksessa henkisten apuv‰lineiden merkitys teknisten rinnalla.

\subsection{Suora ja k‰‰nteinen ongelma}

Inversio-ongelmien eli k‰‰nteis\-ongelmien tutkimuksen tavoitteena on tuottaa molemmanlaisia apuv‰lineit‰ ep‰suoria mittauksia varten.
Tutkin itse k‰‰nteisongelmien matemaattista puolta, joten minulle k‰‰nteisongelmat ovat jokseenkin sama asia kuin ep‰suoran mittaamisen matematiikka.
Keskiˆss‰ ovat usein henkiset apuv‰lineet, ja matematiikan osaamiseksi ei en‰‰ riit‰ kerto- ja jakolaskun hallinta.

Suora ongelma on luonteeltaan seuraavanlainen: kun tiedet‰‰n jonkin esineen ominaisuudet, t‰ytyy p‰‰tell‰, kuinka se k‰ytt‰ytyy.
Jos esimerkiksi kitaran kielest‰ tiedet‰‰n sen pituus, paksuus, tiheys ja j‰nnitys, voidaan p‰‰tell‰, millainen ‰‰ni siit‰ tulee.
T‰ss‰ tapauksessa ‰‰ni koostuu perustaajuudesta, jonka lis‰ksi esiintyy sen monikertoja.
Mainituista tiedoista voi p‰‰tell‰ kaikki mahdolliset taajuudet, joita kieli tuottaa, muttei sit‰, mill‰ voimakkuudella kukin niist‰ esiintyy.

K‰‰nteinen ongelma sen sijaan kysyy saman kysymyksen toiseen suuntaan: kun tiedet‰‰n jonkin esineen k‰ytˆs, t‰ytyy p‰‰tell‰, millaiset sen omainaisuudet ovat.
Kitaraesimerkiss‰ meille siis soitettaisiin kitaran kielt‰, ja kuullusta ‰‰nest‰ pit‰isi p‰‰tell‰ kielen ominaisuudet (pituus, paksuus, tiheys ja j‰nnitys).
T‰m‰ ei kuitenkaan onnistu, sill‰ jos vaikkapa kielen halkaisija puolittuu ja pituus tuplaantuu, pysyy ‰‰ni samanlaisena.
K‰‰nteisongelmalla ei siis olekaan selke‰‰ ratkaisua, mutta kielen ominaisuuksista tiedet‰‰n silti jotain.

Jos kielen nelj‰st‰ ominaisuudesta tiedet‰‰n kolme ja taajuudet ovat selvinneet mittaamalla, voidaan nelj‰s aina p‰‰tell‰.
Jos siis saamme mitata kielen pituuden ja paksuuden ja tied‰mme sen materiaalinkin (voimme katsoa jostain taulukkokirjasta sen tiheyden), saamme selville kielen j‰nnityksen.
T‰ll‰ tavalla muotoiltuna k‰‰nteiseen ongelmaan lˆytyy yksik‰sitteinen ratkaisu.

Mit‰ k‰y, jos kitaran sijasta tutkitaankin rumpua?
Jos rumpukalvon ominaisuudet tiedet‰‰n, voidaan p‰‰tell‰ sen tuottama ‰‰ni (suora ongelma).
Mutta jos kuulemme rummun ‰‰nen ja rumpukalvosta tiedet‰‰n kaikki muu kuin sen muoto, voidaanko muoto p‰‰tell‰?
T‰h‰n kysymykseen ei osaa kukaan antaa t‰ydellist‰ vastausta, mutta jotain sent‰‰n tiedet‰‰n.

Jos ‰‰ni kuulostaa silt‰, ett‰ se tulee t‰ysin pyˆre‰st‰ rummusta, silloin se varmasti tulee pyˆre‰st‰ rummusta.
Pyˆre‰n rummun ‰‰nest‰ ei siis voi erehty‰.
Sen sijaan on olemassa lukuisia esimerkkej‰ erin‰kˆisist‰ kulmikkaista rummuista, joista l‰htee t‰sm‰lleen samanlainen ‰‰ni, joten on myˆs sellaisia rumpuja, joiden ‰‰nest‰ voi erehty‰.

Mit‰ rumpuja ‰‰nen perusteella voi sitten tunnistaa?
T‰t‰ kysymyst‰ tutkitaan kuumeisesti, ja ongelmat ovat nimenomaisesti henkisiss‰ tyˆkaluissa.
T‰m‰ ongelma osataan muotoilla tarkasti matemaattisesti, mutta sen ratkaiseminen on hyvin vaikeaa.
Into t‰h‰n tutkimukseen ei kuitenkaan johdu laajasta rumpuharrastuksesta, vaan samalla matemaattisella mallilla on muitakin sovelluksia.

Esitellyt suorat ja k‰‰nteiset ongelmat edustavat spektraaligeometriaksi kutsuttua alaa.
Nimens‰ mukaisesti se tutkii spektrin (‰‰nen taajuudet) ja geometrian (pituus tai muoto) v‰list‰ suhdetta.
K‰‰nteisongelmia on muunkinlaisia, mutta n‰m‰ k‰yv‰t hyvin esimerkiksi.

Soitinaiheisia esimerkkej‰ muistuttava ‰‰neen perustuva ep‰suora mittaus on monille tuttu.
Vesimelonin viljelij‰t osaavat melonia kumauttamalla ja kuuntelemalla p‰‰tell‰, onko se kyps‰, ja vastaavanlaista p‰‰ttely‰ on osattu k‰ytt‰‰ jo kauan.

K‰‰nteisongelmat ovat tyˆl‰it‰ tutkittavia.
Kuten annetut esimerkit osoittavat, ei ongelmalla v‰ltt‰m‰tt‰ ole yksik‰sitteist‰ ratkaisua (sama k‰ytˆs voi johtua monenlaisista ominaisuuksista) ja ongelma on vaikeampi kuin vastaava suora ongelma.
Varsin kiinnostavaa -- ja sovellusten kannalta hyˆdyllist‰kin -- on jo pelk‰st‰‰n selvitt‰‰, mill‰ k‰‰nteisongelmilla on yksik‰sitteinen ratkaisu eli millaiset ep‰suorat mittausmenetelm‰t voivat toimia edes periaatteessa.

\subsection{Narun punnitus}

Haetaan yhden edellisiin liittyv‰n esimerkin avulla viel‰ lis‰ymm‰rryst‰ ep‰suoraan mittaamiseen.
Otetaan teht‰v‰ksi punnita naru, joka on molemmista p‰ist‰‰n sein‰ss‰ kiinni lˆys‰n pyykkinarun tapaan.
Jos voimme irrottaa narun ja laittaa sen vaa'alle, on ongelma helppo, mutta siit‰ tulee vaikeampi, jos kiell‰mme sen irrottamisen.
T‰llainen rikkomiskielto on luontevampi l‰‰ketieteellisiss‰ sovelluksissa, koska rikottua potilasta ei voi yleens‰ en‰‰ koota, mutta ‰lk‰‰mme antako sen h‰irit‰.

Yksi mahdollisuus on k‰ytt‰‰ narua kitaran kielen‰ ja soveltaa edell‰ esiteltyj‰ oppeja, ja n‰in me teemme.
Sijoitamme ensin l‰helle toista kiinnityskohtaa jonkinlaisen tangon, jonka yli narun annetaan kulkea.
(Jos sopivaa telinett‰ ei lˆydy, pyyd‰mme yst‰v‰‰mme pitelem‰‰n tankoa paikoillaan.)
T‰m‰n tuen ja sein‰n v‰liselle lyhyelle naruosuudelle ripustamme parin kilon punnuksen, jonka massan tied‰mme ja joka on selv‰sti narua painavampi.

Nyt pidempi naruosuus on pingottunut, ja sit‰ voi soittaa kuin kitaran kielt‰.
Vet‰m‰ll‰ narua keskelt‰ v‰h‰n alas ja p‰‰st‰m‰ll‰ liikkeelle saamme narun v‰r‰htelem‰‰n kauniisti.
Mittaamme sekuntikellolla, kauanko narun kest‰‰ tehd‰ kymmenen v‰r‰hdyst‰, ja n‰in saamme selville taajuuden.
Jos koko naru v‰r‰htelee samassa tahdissa ylˆs ja alas, on t‰m‰ taajuus narun perustaajuus.

Narun pituus ja paksuus on helppo mitata.
Lis‰ksi tarvitsemme tiedon j‰nnitysvoimasta.
Sen saamme selville, kun tied‰mme punnuksen massan ja mittaamme, miss‰ kulmassa punnusta kannatteleva naru on.

Nyt olemme aiemmin kuvatussa tilanteessa: tied‰mme narun pituuden, paksuuden ja j‰nnityksen sek‰ sen tuottaman v‰r‰htelytaajuuden.
N‰ist‰ voimme laskea sen tiheyden, ja kertomalla tiheys tilavuudella saadaan selville massa.
T‰m‰ on tietenkin vain pingotetun osan massa, mutta koko narun massa saadaan helposti selville kun katsotaan, kuinka suuri osa koko narusta on pingotettuna.
N‰in on narun massa lopulta saatu selville.

\subsection{Rˆntgen-tomografia}

Er‰‰n k‰‰nteisongelman ratkaisu on helpottanut monen ihmisen el‰m‰‰, ja se ansaitsee maininnan.
Aloitetaan kuitenkin suorasta ongelmasta ja tilanteen kuvailusta.

Rˆntgen-s‰teet ovat hyvin samanlaisia kuin valons‰teet.
Rˆntgen-s‰teily on valoa, jonka v‰ri‰ ihmissilm‰ ei erota -- se on tavallaan liian sinist‰.
T‰llaista valoa voidaan tuottaa ja mitata samaan tapaan kuin lampulla ja kameralla, joten sopivien laitteiden avustamana ihminen kykenee n‰kem‰‰n Rˆntgen-s‰teiden avulla.
Toisin kuin tavallinen valo se kulkee suurimmaksi osaksi ihmisen l‰pi, joten sen avulla voi n‰hd‰ ihmisen sis‰‰n.

Rˆntgen-kuvasta n‰kee helposti murtuneen luun, mutta joitain rakenteita on vaikea n‰hd‰ yhdest‰ kuvasta.
Ja vaikka kuvia olisi montakin, voi olla vaikea n‰hd‰, millainen vaurio esimerkiksi luun sis‰ll‰ on, sill‰ joka suunnasta katsottuna tiell‰ on jotain muuta.

Jos tiedet‰‰n ihmisen tarkka rakenne, eli mit‰ miss‰kin kohdassa kehon sis‰ll‰ on, voidaan t‰st‰ melko helposti p‰‰tell‰, milt‰ mist‰ tahansa suunnasta otettu Rˆntgen-kuva n‰ytt‰‰.
Vastaava k‰‰nteisongelma kysyy, voidaanko kaikista eri suunnista otetuista Rˆntgen-kuvista p‰‰tell‰ ihmisen tarkka kolmiulotteinen rakenne.
T‰m‰ p‰‰ttely on kuin onkin mahdollista, mutta se ei ole helppoa.
P‰‰ttelyyn tarvitaan avuksi tietokone, koska se edellytt‰‰ monimutkaisia laskuja, mutta nykyisilt‰ tietokoneilta t‰m‰ k‰y helposti.
Tietokoneavusteisuuden vuoksi t‰t‰ menetelm‰‰ kutsutaan usein tietokonetomografiaksi (lyhenne CT englanninkielisest‰ nimest‰ computerized tomography).

T‰m‰ ep‰suora mittausmenetelm‰ on eritt‰in hyˆdyllinen.
Joka p‰iv‰ sairaaloissa ymp‰ri maailman l‰‰k‰rit k‰ytt‰v‰t sit‰ n‰hd‰kseen potilaan sis‰‰n.
Ei siis ihme, ett‰ sen keksij‰t Cormack ja Hounsfield saivat Nobelin l‰‰ketieteen palkinnon vuonna 1979.

Oma tutkimukseni liittyy Rˆntgen-tomografian tapaisiin k‰‰nteisongelmiin, mutta sen esittely ei t‰h‰n kirjoitukseen sovi.
Eri alojen k‰‰nteisongelmat kuitenkin liittyv‰t toisiinsa yll‰tt‰vill‰ tavoilla, ja esimerkiksi Rˆntgen-tomografian k‰‰nteisongelma liittyy l‰heisesti edell‰ kuvattuun spektraaligeometriaan.
Juuri t‰m‰ on matematiikan tutkimuksen suola; tuskin kukaan olisi arvannut, ett‰ rummun muodon ja luiden murtumien selvitt‰minen on lopulta hyvin samanlaista!

\end{otherlanguage}

\newpage

\begin{otherlanguage}{latin}
\section{Introductio tironis in problemata inversa}
\label{sec:lat}

Continenter occurrunt rerum condiciones, in quibus aliquid mensurandum est: velocitas autoraedae, capacitas farinae, effectus lucernae aut aliquid aliud.
Quam maxime necesse est instrumenta mensurando apta habere, quo facilius mensurare -- et vivere -- possimus.

In mensurando instrumenta duorum generum habere oportet, technica et mentalia.
Si volumus mensurare velocitatem autoraedae, vehimus ea et mensuramus tempus adhibitum horologio et longitudinem itineris circuitibus canthi numerandis.
Scire volumus et longitudinem unius circuitus, quem funiculo metrico parare possumus.
Sic instrumentis technicis (horologio, numeratro circuituum, funiculo metrico) tres mensurationes factae sunt, sed nulla horum velocitas est.

Deinde adhibenda sunt instrumenta mentalia, artes ratiocinandi et calculandi.
Multiplicatio numeri circuituum eorum longitudine dat longitudinem itineris.
Deinde longitudinem itineris tempore acto dividendo velocitatem mediam cognoscimus.
Sic tandem velocitas mensurata est, in quo maxime usui erant instrumenta et technica et mentalia.

Interdum recte mensurare non licet.
Si medicus densitatem ossis femoris aegri mensurare vult, non licet os e femore removere et deinde massam et capacitatem mensurare de eisque dividendo densitatem calculare.
Salus aegri fovenda limitem ponit rationi mensurandi, et indirecte mensurandum est.
Debet aliquo alio modo densitatem ratiocinari.
Quia ratiocinatio talis haud facilis est, grave est pondus instrumentorum mentalium nec solum technicorum.

\subsection{Problema rectum et inversum}

Scopus investigationis problematum inversorum est instrumenta utriusque generis parare, quae in mensurationibus indirectis adhiberi possint.
Ipse latus mathematicum problematum inversorum scrutor, quare haec problemata artem mathematicam indirecte mensurandi habeo.
In foco sunt instrumenta mentalia, nec, quod ad mathematicen attinet, multiplactio cum divisione satisfacit.

Tale est natura problema rectum: proprietatibus alicuius rei cognitis ratiocinandum est, quo modo illa se gerat.
Si exempli causa notae sunt longitudo, crassitudo, densitas et tensio chordae citharae, vox chorda edita facile calculatur.
Voci enim inest frequentia fundamentalis necnon ea quoque numerorum naturalium multiplicata.
His proprietatibus datis cognosci possunt omnes frequentiae, quas chorda edere potest, sed magnitudo cuiusque frequentiae ignota manet.

Problema inversum autem hanc quaestionem invertit: gestu alicuius rei noto eius proprietates ratiocinandae sunt.
In exemplo nostro citharico vocem audimus et proprietates (longitudo, crassitudo, densitas et tensio) nobis explorandae sunt.
Hoc enim fieri non potest, quia si diametros chordae in dimidium minuitur et longitudo duplicatur, vox eadem manet.
Itaque problema hoc inversum clara solutione caret, sed aliquid de indolibus chordae cognosci potest.

Si tres ex quattuor chordae proprietatibus notae sunt explorataeque frequentiae, quarta semper calculari potest.
Si igitur longitudinem et crassitudinem mensurare licet ac materia nota est (densitas in libro tabulario legi potest), noscitur tensio.
Problema inversum sic formatum solutionem unicam habet.

Quid si citharam in tympanum mutamus?
Si proprietates membranae tympani notae sunt, vox facilis est calculatu (problema rectum).
Sed si vocem tympani audimus et omnes proprietates praeter formam scimus, potestne forma eius explorari?
Nemo responsum perfectum dare quit, sed nonnil iam notum est.

Si vox tympani talis est, ut eam membrana rotunda dare possit, certe membrana rotunda data est.
Vox tympani rotundi numquam nos fallit.
Cum autem multae formae angulosae, quae eandem vocem dent, inventae sint, possumus etiam errare de aliquibus formis tympanorum.

Quales tympanos possumus eorum vocibus auditis recognoscere?
Haec quaestio acerrime investigatur, et problemata sunt in instrumentis mentalibus.
Haec quaestio mathematice exacte formulari potest, sed eam solvere difficillimum est.
Studium huius quaestionis non ab oblectationibus musicalibus oritur, sed idem exemplar mathematicum alias quoque usui est.

Problemata recta et inversa supra exposita pertinent ad aream mathematices, quae geometria spectralis appellatur.
Secundum nomen suum nexum inter spectrum (frequentiae vocis) et geometriam (longitudo sive forma) scrutatur.
Alia quoque sunt problemata inversa, sed haec bonum exemplum praebent.

Mensuratio indirecta exemplis musicalibus similis multis nota est.
Cultores melonum melonem quatiendo et audiendo ratiocinari possunt, an maturus sit, et talis ratiocinatio iam diu adhibita est.

Problemata inversa investiganda laboriosa sunt.
Ut exempla demonstrant, interdum solutio unica problemati deest (variae proprietates eosdem gestus efficere possunt) et problema est difficilius quam versio eius recta.
Interest -- et usui est -- investigare, quae problemata inversa solutionem unicam habeant sive quae methodi indirecte mensurandi omnino adhiberi possint.

\subsection{Restis pensanda}

Conemur uno exemplo mensurationes indirectas melius intelligere.
Metam nobis ponamus restem pensandam, quae utroque limite parieti affixa est sicut restis lintearia.
Si restem abripere licet et ad libram ponere, problema facile est, sed difficilius fit, si eam non licet solvere.
Prohibitio talis naturalior est in arte medicina, cum homo aegrotus fractus restitui non possit, sed ne hoc distrahamur.

Possumus restem sicut chordam citharae adhibere et meminisse supra scripta, et sic facimus.
Ponimus prope alterum locum affixionis axem alicuius generis, super quam restis eat.
(Si nil est, in quo axem affigeremus, rogamus amicum, ut axem teneat.)
In parte breviore restis suspendimus pondus duorum fere chiliogrammatum, cuius massa est multo maior quam restis et nobis nota.

Nunc pars longior restis tenta est, et ea canere possumus sicut chorda citharae.
Si mediam restem paullo deorsum trahimus deinque solvimus, restis pulchre oscillat.
Horologio mensuramus quanto tempore decies oscillet et ita frequentiam cognoscimus.
Tota reste eodem rhythmo oscillante scimus hanc frequentiam esse frequentiam fundamentalem.

Longitudo et crassitudo sunt faciles mensuratu.
Volumus scire etiam tensionem.
Eam cognoscimus, quia massa ponderis nota est et angulus partis restis, de qua pendet, facile mensuratur.

Nunc sumus in rerum condicione supra descripta: scimus longitudinem, crassitudinem, tensionem et frequentiam oscillationis.
Ex his calculare possumus densitatem, et illam capacitate multiplicando massam exploramus.
Haec est solum massa partis restis tentae, sed massam totius restis facile cognoscimus videndo, qualis pars restis tenta sit.
Sic demum est massa restis cognita.

\subsection{Tomographia Rˆntgeniana}

Solutio cuiusdam problematis inversi vitam multorum faciliorem reddit, et mentionis dignum est.
Incipiamus enim cum problemate directo et descriptione rerum condicionum.

Radii Rˆntgeniani persimiles sunt radiis lucis.
Radiatio Rˆntgeniana est lux, cuius colorem oculus humanus non capit -- est, ut ita dicam, nimis caeruleus.
Lux talis parari et mensurari potest sicut lux visibilis lampade et machina photographica, quare instrumentis aptis ornatus homo radiis Rˆntgenianis videre potest sicut solitis radiis lucis.
Dissimiliter luci hodiernae haec lux magnopere per corpus humanum it, quod efficit, ut interiora hominis videre liceat.

In imagine Rˆntgeniana os fractum facile videtur, sed non omnes structurae tam facile in una imagine percipiuntur.
Etsi multae sint imagines, est defectus in medulla ossis plerumque difficilis perceptu, cum ossibus aliisque structuris a quaque parte umbretur.

Si structura hominis accurate nota est, id est, quid ubique sit, facile est ratiocinari, quales imagines Rˆntgenianae a quaque parte factae videantur.
Problema inversum ad hoc attinens rogat, possitne structura hominis tridimensionalis ex imaginibus Rˆntgenianis ab omnibus partibus factis inveniri.
Accidit, ut vero inveniri possit, sed haud facile est.
In ratiocinando auxilio computatri opus est, quia computationes sunt complicatae, sed computatris hodiernis hoc facillime efficitur.
Propter usum computatri haec methodus plerumque tomographia computatralis appellatur (breviter CT ex verbis Anglicis computerized tomography).

Haec methodus indirecta mensurandi utilissima est.
Cotidie in nosocomiis medici ea utuntur, ut intra aegros videant.
Itaque haud mirum est eius inventores Cormack et Hounsfield praemium Nobelianum in medicina anno 1979 tulisse.

Investigationes meae ad problemata inversa tomographiae Rˆntgenianae similia attinent, sed expositio earum hic non apta est.
Attamen diversa problemata inversa nexus improvisos inter se habent, et exempli causa problema inversum tomographiae Rˆntgenianae nexum intimum habet cum geometria spectrali supra exposita.
Hoc est sal mathematices; haud quisquam coniectaverit ratiocinationem formae tympani et fracturae ossis simillimam esse!

\end{otherlanguage}
\end{sloppypar}
\end{appendices}

\end{document}